\documentclass[12pt]{article}
 \usepackage{amssymb}
    \usepackage{latexsym}

\title{{\bf Twisted modules for vertex operator algebras}}
    \author{Benjamin Doyon\footnote{
Department of Mathematical Sciences, Durham University, U.K.}{\ }\footnote{Work done mainly while at Rudolf Peierls Centre for Theoretical Physics, Oxford University, U.K, under EPSRC Postdoctoral Fellowship, grant GR/S91086/01.}}
    \date{}
    
\begin{document}

    \bibliographystyle{alpha}
    \maketitle
\begin{abstract}
This contribution is mainly based on joint papers with Lepowsky and Milas, and some parts of these papers are reproduced here. These papers further extended works by Lepowsky and by Milas. Following our joint papers, I explain the general principles of twisted modules for vertex operator algebras in their powerful formulation using formal series, and derive general relations satisfied by twisted and untwisted vertex operators. Using these, I prove new ``equivalence'' and ``construction'' theorems, identifying a set of sufficient conditions in order to have a twisted module for a vertex operator algebra, and a simple way of constructing the twisted vertex operator map. This essentially combines our general relations for twisted modules with ideas of Li (1996), who had obtained similar construction theorems using different relations. Then, I show how to apply these theorems in order to construct twisted modules for the Heisenberg vertex operator algebra. I obtain in a new way the explicit twisted vertex operator map, and in particular give a new derivation and expression for the formal operator $\Delta_x$ constructed some time ago by Frenkel, Lepowsky and Meurman. Finally, I reproduce parts of our joint papers. I use the untwisted relations in the Heisenberg vertex operator algebra in order to understand properties of a certain central extension of a Lie algebra of differential operators on the circle: the connection between the structure of the central term in Lie brackets and the Riemann Zeta function at negative integers. I then use the twisted relations in order to construct in a simple way a family of representations for this algebra based on twisted modules for the Heisenberg vertex operator algebra. As a simple consequence of the twisted relations, the construction involves the Bernoulli polynomials at rational values in a fundamental way.
\end{abstract}

    \input amssym.def
    \input amssym

    \newtheorem{theorem}{Theorem}[section]
    \newtheorem{remark}{Remark}
    \newtheorem{rema}{Remark}[section]
    \newtheorem{propo}[rema]{Proposition}
    \newtheorem{theo}[rema]{Theorem}
    \newtheorem{defi}[rema]{Definition}
    \newtheorem{lemma}[rema]{Lemma}
    \newtheorem{corol}[rema]{Corollary}

\renewcommand{\theequation}{\thesection.\arabic{equation}}
\renewcommand{\theremark}{\thesection.\arabic{remark}}
\renewcommand{\thetheo}{\thesection.\arabic{remark}}
\renewcommand{\thedefi}{\thesection.\arabic{remark}}
\renewcommand{\thecorol}{\thesection.\arabic{remark}}
\setcounter{section}{0} \setcounter{equation}{0}
\setcounter{remark}{0}

\def\sect#1{\section{#1}\setcounter{equation}{0}\setcounter{rema}{0}}
\def\ssect#1{\subsection{#1}}
\def\sssect#1{\subsubsection{#1}}

\newcommand{\bc}{\begin{center}}
\newcommand{\ec}{\end{center}}
\newcommand{\ba}{\begin{array}}
\newcommand{\ea}{\end{array}}
\newcommand{\nn}{\nonumber \\}
\newcommand{\beq}{\begin{equation}}
\def\sop{x_1^{1/p} \rightarrow \omega_p^s (x_2+x_0)^{1/p}}
\def\rp{x_1^{1/p} \rightarrow \omega_p^r (x_2+x_0)^{1/p}}
\def\xp{x^{1/p} \rightarrow \omega_p^s x^{1/p}}
\newcommand{\eeq}{\end{equation}}
\newcommand{\beqa}{\begin{eqnarray}}
\newcommand{\eeqa}{\end{eqnarray}}
\newcommand{\no}{\nonumber}
\def\bi{\begin{itemize}}
\def\ei{\end{itemize}}
\def\mato#1{\left(\ba{#1}} 
\def\matf{\ea\right)}

\def\eq#1{(\ref{#1})}
\def\lab#1{\label{#1}}

\def\d{\partial}
\def\dv#1#2{\frac{\delta #1}{\delta #2}}
\def\frc#1#2{\frac{#1}{#2}}
\def\b#1{\bar{#1}}
\def\t#1{\tilde{#1}}
\def\h#1{\hat{#1}}
\def\lt#1{\left#1}
\def\rt#1{\right#1}
\def\la{\langle}
\def\ra{\rangle}
\def\<{\langle}
\def\>{\rangle}
\def\:{\mbox{\tiny ${\bullet\atop\bullet}$}}

\def\F{{\Bbb F}}
\def\Z{{\Bbb Z}}
\def\C{{\Bbb C}}
\def\R{{\Bbb R}}
\def\N{{\Bbb N}}
\def\D{{\cal D}}
\def\a#1{{\goth{#1}}}
\def\End{{\rm End\,}}
\def\Hom{{\rm Hom}}
\def\Der{{\rm Der}}
\def\mod{{\rm \,mod\,}}
\def\dim{{\rm dim\,}}
\def\Res{{\rm Res}}
\def\Tr{{\rm Tr}}
\def\proof{{\em Proof: }}
\newcommand{\halmos}{\rule{1ex}{1.4ex}}
\newcommand{\eproof}{\hspace*{\fill}\mbox{$\halmos$}}

\def\for#1{\quad\mbox{#1 }}
\def\com#1{\quad(#1)}
\def\om{\omega_p}

\newcommand{\nordplus}{\mbox{\scriptsize ${+ \atop +}$}}
\newcommand{\nordbullet}{\mbox{\tiny ${\bullet\atop\bullet}$}}

\def\aa{\mbox{\scriptsize ${\ast\atop\ast}$}}
\def\pp{\nordplus}
\def\rr{\nordbullet}
\def\xx{\mbox{\scriptsize ${\times\atop\times}$}}

\def\glim#1{\widehat{\lim_{#1}}}

\def\dpf#1#2{\delta\left(\frac{#1}{#2}\rt)}
\def\pf#1#2{\left(\frac{#1}{#2}\rt)}
\def\id{1_\a{h}}
\sect{Introduction}

This contribution is mainly based on, and partly reproduces, the recent papers by the
present author, Lepowsky and Milas \cite{DLMi1}, \cite{DLMi2}.
These works were a continuation of a series of papers of Lepowsky
and Milas \cite{L3}, \cite{L4}, \cite{M1}--\cite{M3}, stimulated
by work of Bloch \cite{Bl}.

In those papers, we used the general theory of vertex operator
algebras to study central extensions of classical Lie algebras and
superalgebras of differential operators on the circle in
connection with values of $\zeta$--functions at the negative
integers, and with the Bernoulli polynomials at rational values.
Parts of the present contribution recall the main results of
\cite{DLMi1,DLMi2}: Using general principles of the theory of
vertex operator algebras and their twisted modules, we obtain a
bosonic, twisted construction of a certain central extension of a
Lie algebra of differential operators on the circle, for an
arbitrary twisting automorphism. The construction involves the
Bernoulli polynomials in a fundamental way.  This is explained
through results in the general theory of vertex operator algebras,
including an identity discovered in \cite{DLMi1,DLMi2} which was called 
``modified weak associativity'', and which is a consequence of the twisted Jacobi identity.

More precisely, we combine and extend methods {}from \cite{L3}, \cite{L4},
\cite{M1}--\cite{M3}, \cite{FLM1}, \cite{FLM2} and \cite{DL2}. In
those earlier papers, vertex operator techniques were used to
analyze untwisted actions of the Lie algebra
$\hat{\mathcal{D}}^+$, studied in \cite{Bl}, on a module for a
Heisenberg Lie algebra of a certain standard type, based on a
finite-dimensional vector space equipped with a nondegenerate
symmetric bilinear form. Now consider an arbitrary isometry $\nu$
of period say $p$, that is, with $\nu^p={1}$. Then, it was shown in
\cite{DLMi1,DLMi2} that the corresponding $\nu$--twisted modules carry
an action of the Lie algebra $\hat{\mathcal{D}}^+$ in terms of
twisted vertex operators, parametrized by certain quadratic
vectors in the untwisted module. This extends a result
{}from \cite{FLM1}, \cite{FLM2}, \cite{DL2} where actions of the
Virasoro algebra were constructed using twisted vertex operators.

Still following \cite{DLMi1,DLMi2}, we explicitly compute certain ``correction'' terms for
the generators of the ``Cartan subalgebra'' of
$\hat{\mathcal{D}}^+$ that naturally appear in any twisted
construction. These correction terms are expressed in terms of
special values of certain Bernoulli polynomials. They can in
principle be generated, in the theory of vertex operator algebras,
by the formal operator $e^{\Delta_x}$ \cite{FLM1}, \cite{FLM2},
\cite{DL2} involved in the construction of a twisted action for a
certain type of vertex operator algebra, the Heisenberg vertex
operator algebra. We generate those correction terms in an easier
way, using the modified weak associativity relation.

Then, the present contribution extends the works
\cite{DLMi1,DLMi2} described above by providing a detailed analysis of the
modified weak associativity relation. We state and prove a new
theorem (Theorem \ref{theoequil}) about the equivalence of modified weak associativity and
weak commutativity with the twisted Jacobi identity, and a new
``construction'' theorem (Theorem \ref{theoconstr1}), where we identify a set
of sufficient conditions in order to have a twisted module for a
vertex operator algebra, and a simple way of constructing the
twisted vertex operator map. The latter theorem essentially combines modified weak associativity with ideas of Li \cite{Li1,Li2}, where similar construction theorems were proven using different general relations of vertex operator algebras and twisted modules -- there may be a ``direct'' path from Li's construction theorems to ours, but we haven't investigated this. The use of modified weak associativity seems to have certain advantages in the twisted case. As an illustration, we give a new proof that the $\nu$-twisted Heisenberg Lie algebra modules mentionned above are also twisted modules for the Heisenberg vertex operator algebra. Using our theorems, we explicitly construct the twisted vertex operator map (Theorem \ref{theotvom}). This gives a new and relatively simple derivation and expression for this map, and in particular for the formal operator $\Delta_x$ mentioned above. A consequence of this is that one minor technical assumption that had to be made in \cite{DLMi1,DLMi2}, about the action of the automorphism $\nu$, can be taken away.

We should mention that in \cite{KR} Kac and Radul established a relationship between the
Lie algebra of differential operators on the circle and the Lie
algebra $\widehat{\frak{gl}}(\infty)$; for further work in this
direction, see \cite{AFMO}, \cite{KWY}. Our methods and motivation
for studying Lie algebras of differential operators, based on
vertex operator algebras, are new and very different, so we do not
pursue their direction.

Although we will present many of the main results of \cite{DLMi2} with some of
the proofs, we refer the reader to this paper for a more extensive
discussion of those results.

{\em Acknowledgments.} The author is grateful to J. Lepowsky and A. Milas for discussions
and comments on the manuscript.

\sect{Vertex operator algebras, untwisted modules and twisted
modules} \label{sectVOA}

In this section, we recall the definition of vertex operator
algebras, (untwisted) modules and twisted modules. For the basic
theory of vertex operator algebras and modules, we will use the
viewpoint of \cite{LL}.

In the theory of vertex operator algebras, formal calculus plays a
fundamental role.  Here we recall some basic elements of formal
calculus (cf. \cite{LL}).  Formal calculus is the calculus of
formal doubly--infinite series of formal variables, denoted below
by $x$, $y$, and by $x_1,\,x_2,\ldots$, $y_1,y_2,\ldots$.  The
central object of formal calculus is the formal delta--function
\[
    \delta(x) = \sum_{n\in\Z} x^n
\]
which has the property
\[
    \delta\lt(\frc{x_1}{x_2}\rt) f(x_1) = \delta\lt(\frc{x_1}{x_2}\rt)
    f(x_2)
\]
for any formal series $f(x_1)$. The formal delta--function enjoys
many other properties, two of which are:
\beq\label{2delta}
    x_2^{-1}\delta\lt(\frc{x_1-x_0}{x_2}\rt) = x_1^{-1}
    \delta\lt(\frc{x_2+x_0}{x_1}\rt)
\eeq
and
\beq\label{3delta}
    x_0^{-1}\delta\lt(\frc{x_1-x_2}{x_0}\rt) + x_0^{-1}
    \delta\lt(\frc{x_2-x_1}{-x_0}\rt) =
    x_2^{-1}\delta\lt(\frc{x_1-x_0}{x_2}\rt).
\eeq
In these equations, binomial expressions of the type
$(x_1-x_2)^n,\,n\in\Z$ appear. Their meaning as formal series in
$x_1$ and $x_2$, as well as the meaning of powers of more
complicated formal series, is summarized in the ``binomial
expansion convention'' -- the notational device according to which
binomial expressions are understood to be expanded in nonnegative
integral powers of the second variable. When more elements of
formal calculus are needed below, we shall recall them.

\ssect{Vertex operator algebras and untwisted modules}

We recall {}from \cite{FLM2} the definition of the notion of
vertex operator algebra, a variant of Borcherds' notion \cite{Bo}
of vertex algebra:
\begin{defi}\label{VOA}
A {\bf vertex operator algebra} $(V,Y,{\bf 1},\omega)$, or $V$ for
short, is a $\mathbb{Z}$--graded vector space
\[
    V=\coprod_{n \in \mathbb{Z}} V_{(n)};
    \ \mbox{\rm for}\ v\in V_{(n)},\;\mbox{\rm wt}\ v = n,
\]
such that
\beqa
     &&  V_{(n)} = 0 \;\; \mbox{ for }\; n \mbox{ sufficiently
    negative,} \nn &&  \mbox{\rm dim }V_{(n)}<\infty\;\;\mbox{ for }\;
    n \in {\mathbb Z} , \no
\eeqa
equipped with a linear map $Y(\cdot,x)$:
\begin{eqnarray}
    Y(\cdot,x)\,: \ V&\to&(\mbox{\rm End}\; V)[[x, x^{-1}]]\nonumber \\
    v&\mapsto& Y(v, x)={\displaystyle \sum_{n\in{\mathbb Z}}}v_{n} x^{-n-1}
    \,,\;\;v_{n}\in \mbox{\rm End} \;V, \label{voalg}
\end{eqnarray}
where $Y(v,x)$ is called the {\em vertex operator} associated with
$v$, and two particular vectors, ${\bf 1},\,\omega\in V$, called
respectively the {\em vacuum vector} and the {\em conformal
vector}, with the following properties:\\
{\em truncation condition:} For every $v,w \in V$
\beq
    v_n w=0
\eeq
for $n\in\Bbb Z$ sufficiently large;\\
{\em vacuum property:}
\beq
    Y({\bf 1},x) = 1_V \for{($1_V$ is the identity on $V$);}
\eeq
{\em creation property:}
\beq
    Y(v,x){\bf 1} \in V[[x]] \for{and} \lim_{x\to0} Y(v,x){\bf 1} = v \;;
\eeq
{\em Virasoro algebra conditions:} Let
\begin{equation}
    L(n)=\omega _{n+1}\;\; \mbox{\rm for} \;n\in{\mathbb Z},
        \;\;{\rm i.e.},\;\;
        Y(\omega, x)=\sum_{n\in{\mathbb Z}}L(n)x^{-n-2} \;.
\end{equation}
Then
\[
    [L(m),L(n)]=
    (m-n)L(m+n)+c_V\frc{m^{3}-m}{12}\,\delta_{n+m,0}\,1_V
\]
for $m, n \in {\mathbb Z}$, where $c_V\in {\mathbb C}$ is the
central charge (also called ``rank'' of $V$),
\[
    L(0) v=({\rm wt}\ v) v
\]
for every homogeneous element $v$, and we have the {\em
$L(-1)$--derivative property:}
\beq\label{Lm1der}
    Y(L(-1)u,x)=\frac{d}{dx}Y(u,x) \;;
\eeq
{\em Jacobi identity:}
\beqa\label{jacobi}
    && x_0^{-1}\dpf{x_1-x_2}{x_0} Y(u,x_1) Y(v,x_2)
        - x_0^{-1}\dpf{x_2-x_1}{-x_0} Y(v,x_2) Y(u,x_1) \no\\
    &&  \qquad = x_2^{-1} \delta\lt(\frc{x_1-x_0}{x_2}\rt)
            Y(Y(u,x_0)v,x_2)\;.
\eeqa
\end{defi}

An important property of vertex operators is skew--symmetry, which
is an easy consequence of the Jacobi identity (cf. \cite{FHL}):
\beq \label{skewsymm}
    Y(u,x)v = e^{x L(-1)} Y(v,-x)u.
\eeq
Another easy consequence of the Jacobi identity is the
$L(-1)$--bracket formula:
\beq\label{Lm1bracket}
    [L(-1),Y(u,x)] = Y(L(-1)u,x).
\eeq

Fix now a vertex operator algebra $(V,Y,{\bf 1},\omega)$, with
central charge $c_V$.
\begin{defi}\label{VOAmodule}
A (${\Bbb Q}$--graded) {\bf module} $W$ for the vertex operator
algebra $V$ (or $V$--{\bf module}) is a ${\Bbb Q}$--graded vector
space,
$$
    W=\coprod_{n\in {\Bbb Q}}W_{(n)}; \for{for} v\in
    W_{(n)},\;\mbox{\rm wt}\ v = n,
$$
such that \beqa &&  W_{(n)} = 0 \for{for} n \mbox{ sufficiently
negative,} \nn &&  \dim W_{(n)}<\infty\for{for} n \in {\Bbb Q} ,
\no \eeqa equipped with a linear map
\begin{eqnarray}
    Y_W(\cdot,x)\,: \ V&\to&(\mbox{\rm End}\; W)[[x, x^{-1}]]\nonumber \\
    v&\mapsto& Y_W(v, x)={\displaystyle \sum_{n\in{\Bbb Z}}}v_{n}^W x^{-n-1}
    \,,\;\; v_{n}^W\in \mbox{\rm End} \; W,
\label{vo}
\end{eqnarray}
where $Y_W(v,x)$ is still called the {\em vertex operator}
associated with $v$, such that the following conditions hold:\\
{\em truncation condition:} For every $v \in V$ and $w \in W$
\beq\label{truncation}
    v_n^W w=0
\eeq
for $n\in\Bbb Z$ sufficiently large;\\
{\em vacuum property:}
\begin{equation} \label{vacuum}
    Y_W({\bf 1}, x)=1_{W};
\end{equation}
{\em Virasoro algebra conditions:} Let
\[
    L_W(n)=\omega _{n+1}^W\;\; \mbox{\rm for} \;n\in{\mathbb Z},
    \;\;{\rm i.e.},\;\;
    Y_W(\omega,x)=\sum_{n \in \mathbb{Z}} L_W (n)x^{-n-2}.
\]
We have
\beqa
    [L_W(m),L_W(n)]&=&(m-n)L_W(m+n)+
        c_V\frac{m^3-m}{12}\,\delta_{m+n,0}\, 1_W,\nn
    L_W(0) v&=&({\rm wt}\ v) v \no
\eeqa
for every homogeneous element $v\in W$, and
\beq\label{Lm1derm}
        Y_W(L(-1)u,x)=\frac{d}{dx}Y_W(u,x) \;;
\eeq
{\em Jacobi identity:}
\beqa\label{jacobim}
    &&  x_0^{-1}\dpf{x_1-x_2}{x_0} Y_W(u,x_1) Y_W(v,x_2)
    - x_0^{-1}\dpf{x_2-x_1}{-x_0} Y_W(v,x_2) Y_W(u,x_1) \no\\
    &&  \qquad = x_2^{-1} \delta\lt(\frc{x_1-x_0}{x_2}\rt)
        Y_W(Y(u,x_0)v,x_2).
\eeqa
\end{defi}

{}From the Jacobi identity (\ref{jacobim}), one can derive the
weak commutativity and weak associativity relations, respectively:
\beqa
    (x_1-x_2)^{k(u,v)}Y_W(u,x_1) Y_W(v,x_2) &=& (x_1-x_2)^{k(u,v)}
    Y_W(v,x_2) Y_W(u,x_1) \no \\
    \label{wcom} \\
    (x_0+x_2)^{l(u,w)} Y_W(u,x_0+x_2) Y_W(v,x_2)w &=& (x_0+x_2)^{l(u,w)}
    Y_W(Y(u,x_0)v,x_2)w,
    \no \\ \label{wass}
\eeqa where $u,\,v\,\in V$ and $w\in W$, valid for large enough
nonnegative integers $k(u,v)$ and $l(u,w)$, their minimum value depending
respectively on $u,\,v$ and on $u,\,w$. For definiteness, we will
pick the integers $k(u,v)$ and $l(u,w)$ to be the smallest
nonnegative integers for which the relations above are valid.

\ssect{Twisted modules for vertex operator algebras}

The notion of twisted module for a vertex operator algebra was
formalized in \cite{FFR} and \cite{D} (see also the geometric
formulation in \cite{FrS}; see also \cite{DLM}), summarizing the
basic properties of the actions of twisted vertex operators
discovered in \cite{FLM1}, \cite{FLM2} and \cite{L2}; the main
nontrivial axiom in this notion is the twisted Jacobi identity of
\cite{FLM2} (and \cite{L2}); cf. \cite{FLM1}.

A critical ingredient in formal calculus needed in the theory of
twisted modules is the appearance of fractional powers of formal
variables, like $x^{1/p},\, p\in \Z_+$ (the positive integers).
For the purpose of formal calculus, the object $x^{1/p}$ is to be
treated as a new formal variable whose $p$--th power is $x$.  The
binomial expansion convention is applied as stated at the
beginning of Section \ref{sectVOA} to binomials of the type
$(x_1+x_2)^{1/p}$. {}From a geometrical point of view, these rules
correspond to choosing a branch in the ``orbifold structure''
described (locally) by the twisted vertex operator algebra module.

We now fix a positive integer $p$ and a primitive $p$--th root of
unity
\beq\label{primitiveroot}
    \om \in \mathbb{C}.
\eeq
We record here two important properties of the formal
delta--function involving fractional powers of formal variables:
\beq\label{deltaid1}
    \delta(x) = \frc1p\sum_{r=0}^{p-1}\delta(\omega^{r}_p
    x^{1/p})
\eeq
and
\beq\label{deltaid2}
    x_2^{-1} \delta\lt(\om^r\lt(\frc{x_1-x_0}{x_2}\rt)^{1/p}\rt) =
    x_1^{-1}
    \delta\lt(\om^{-r}\lt(\frc{x_2+x_0}{x_1}\rt)^{1/p}\rt).
\eeq
The latter formula can be found (in a slightly different form) in \cite{Li2}.
For the sake of completeness, we present here a proof.

\proof The coefficient of $x_0^0$ in equation (\ref{deltaid2})
is immediate. Consider some formal series $f(x) = \sum_{n\in\C} f_n x^n,\; f_n\in\C$. From the formula
\[
	(-1)^k  (\d/\d x_1)^k (x_1^s x_2^{-s-1}) = (\d/\d x_2)^k (x_1^{s-k} x_2^{-s-1+k})
\]
for any $s\in\C$ and $k$ a nonnegative integer, we find that
\beq
	(-1)^k \lt(\frc{\d}{\d x_1}\rt)^k \lt(x_2^{-1} f(x_1/x_2)\rt) =  \lt(\frc{\d}{\d x_2}\rt)^k \lt(x_1^{-1} (x_1/x_2)^{1-k} f(x_1/x_2)\rt).
\eeq
With $f(x) = \delta\lt(\om^r x^{1/p}\rt)$, we use the formal delta-function property to get
\[
	(x_1/x_2)^{1-k}\delta\lt(\om^r (x_1/x_2)^{1/p}\rt)= \delta\lt(\om^r (x_1/x_2)^{1/p}\rt)
\]
and thus
\beq
	(-1)^k \lt(\frc{\d}{\d x_1}\rt)^k \lt(x_2^{-1} \delta\lt(\om^r(x_1/x_2)^{1/p}\rt)\rt) =  \lt(\frc{\d}{\d x_2}\rt)^k \lt(x_1^{-1} \delta\lt(\om^{-r}(x_2/x_1)^{1/p}\rt)\rt).
\eeq
Summing over nonnegative integers $k$ with the coefficients $x_0^k/k!$ on both sides, we obtain (\ref{deltaid2}). \eproof

Recall the vertex operator algebra $(V, Y, {\bf 1}, \omega)$ with
central charge $c_V$ of the previous subsection.  Fix an
automorphism $\nu$ of period $p$ of the vertex operator algebra
$V$, that is, a linear automorphism of the vector space $V$
preserving $\omega$ and ${\bf 1}$ such that
\beq
    \nu Y(v,x)\nu^{-1} = Y(\nu v,x) \ \mbox{for}\ v\in V,
\eeq
and
\beq
    \nu^p=1_V.
\eeq

\begin{defi}\label{VOAmodulet}
A (${\Bbb Q}$-graded) $\nu$-{\bf twisted} $V$-{\bf module} $M$ is
a ${\Bbb Q}$-graded vector space,
$$
    M=\coprod_{n\in {\Bbb Q}}M_{(n)}; \ \mbox{\rm for}\ v\in
    M_{(n)},\;\mbox{\rm wt}\ v = n,
$$
such that
\beqa
    &&  M_{(n)} = 0 \;\; \mbox{ for }\; n \mbox{ sufficiently negative,} \nn
    &&  \mbox{\rm dim }M_{(n)}<\infty\;\;\mbox{ for }\; n \in {\Bbb Q} , \no
\eeqa
equipped with a linear map
\begin{eqnarray}
    Y_M(\cdot,x)\,: \ V&\to&(\mbox{\rm End}\; M)[[x^{1/p},
    x^{-1/p}]]\nonumber \\
    v&\mapsto& Y_M(v, x)={\displaystyle \sum_{n\in{\frc1p\Bbb
    Z}}}v_{n}^\nu x^{-n-1}
    \,,\;\; v_{n}^\nu\in \mbox{\rm End} \; M,
\label{tvo}
\end{eqnarray}
where $Y_M(v,x)$ is called the {\em twisted vertex operator}
associated with $v$, such that the following conditions hold:\\
{\em truncation condition:} For every $v \in V$ and $w \in M$
\beq\label{truncationt}
    v_n^\nu w=0
\eeq
for $n\in\frc1p\Bbb Z$ sufficiently large;\\
{\em vacuum property:}
\begin{equation} \label{vacuumt}
    Y_M({\bf 1}, x)=1_{M};
\end{equation}
{\em Virasoro algebra conditions:} Let
\[
    L_M(n)=\omega _{n+1}^\nu\;\; \mbox{\rm for} \;n\in{\mathbb Z},
    \;\;{\rm i.e.},\;\;
    Y_M(\omega,x)=\sum_{n \in \mathbb{Z}} L_M (n)x^{-n-2}.
\]
We have
\beqa
    [L_M(m),L_M(n)]&=&(m-n)L_M(m+n)+
        c_V\frac{m^3-m}{12}\,\delta_{m+n,0}\, 1_M,\nn
    L_M(0) v&=&({\rm wt}\ v) v
\label{L0weightmt}
\eeqa
for every homogeneous element $v$, and
\beq\label{Lm1dermt}
    Y_M(L(-1)u,x)=\frac{d}{dx}Y_M(u,x)\;;
\eeq
{\em Jacobi identity:}
\beqa\label{jacobitm}
    && x_0^{-1}\dpf{x_1-x_2}{x_0} Y_M(u,x_1) Y_M(v,x_2)
    - x_0^{-1}\dpf{x_2-x_1}{-x_0} Y_M(v,x_2) Y_M(u,x_1) \no\\
    &&  \qquad = \frc1p x_2^{-1} \sum_{r=0}^{p-1}
    \delta\lt(\om^r\lt(\frc{x_1-x_0}{x_2}\rt)^{1/p}\rt)
    Y_M(Y(\nu^r u,x_0)v,x_2).
\eeqa
\end{defi}

Note that when restricted to the fixed--point subalgebra $\{u \in
V\,|\, \nu u = u\}$, a twisted module becomes a {\em true} module:
the twisted Jacobi identity (\ref{jacobitm}) reduces to the
untwisted one (\ref{jacobim}), by (\ref{deltaid1}). This will
enable us to construct natural representations of a certain
infinite-dimensional algebra $\h{\D}^+$ (see Section \ref{sect-D})
on suitable twisted modules.

\sect{Heisenberg vertex operator algebra and its twisted modules}
\label{sect-Heisenberg}

It is appropriate at this point to make these definitions more
substantial by giving a simple but important example of a vertex
operator algebra, and of some of its twisted modules.

\ssect{Heisenberg vertex operator algebra}

Following \cite{FLM2}, let $\a{h}$ be a finite-dimensional abelian
Lie algebra (over $\C$) of dimension $d$ on which there is a
nondegenerate symmetric bilinear form $\<\cdot,\cdot\>$.  Let
$\nu$ be an isometry of $\a{h}$ of period $p>0$:
$$
    \<\nu\alpha,\nu\beta\> = \<\alpha,\beta\> ,\quad
    \nu^p\alpha = \alpha
$$
for all $\alpha,\beta\in\a{h}$. Consider the affine Lie algebra
$\hat{\a{h}}$,
\[
    \h{\a{h}} = \coprod_{n\in\Z} \a{h} \otimes t^n \oplus \C C,
\]
with the commutation relations
\begin{eqnarray}
    [\alpha\otimes t^m,\beta\otimes t^n] &=& \<\alpha,\beta\>
    m\delta_{m+n,0} C \com{\alpha,\beta \in\a{h},\; m,n\in\Z} \no\\ {}
    [C,\h{\a{h}}] &=& 0. \no
\end{eqnarray}
Set
\[
    \h{\a{h}}^+ = \coprod_{n>0} \a{h}\otimes t^n ,\qquad
    \h{\a{h}}^- = \coprod_{n<0} \a{h}\otimes t^n.
\]
The subalgebra
\[
    \h{\a{h}}^+ \oplus \h{\a{h}}^- \oplus \C C
\]
is a Heisenberg Lie algebra. Form the induced (level--one)
$\h{\a{h}}$--module
\[
    S = \mathcal{U}(\hat{\a{h}}) \otimes_{\mathcal{U}\lt(
    \hat{\a{h}}^+ \oplus \a{h} \oplus \C C\rt)} \C
        \simeq S(\hat{\a{h}}^-) \for{(linearly),}
\]
where $\hat{\a{h}}^+\oplus \a{h}$ acts trivially on $\C$ and $C$
acts as 1; $\mathcal{U}(\cdot)$ denotes universal enveloping
algebra and $S(\cdot)$ denotes the symmetric algebra.  Then $S$ is
irreducible under the Heisenberg algebra $\h{\a{h}}^+ \oplus
\h{\a{h}}^- \oplus \C C$.  We will use the notation
$\alpha(n)\;(\alpha\in\a{h},\, n\in\Z)$ for the action of
$\alpha\otimes t^n \in \hat{\a{h}}$ on $S$.

The induced $\hat{\a{h}}$-module $S$ carries a natural structure
of vertex operator algebra. This structure is constructed as
follows (cf. \cite{FLM2}). First, one identifies the vacuum vector
as the element 1 in $S$: ${\bf 1}=1$. Consider the following
formal series acting on $S$:
$$
    \alpha(x)=\sum_{n \in \Z} \alpha(n)x^{-n-1}\com{\alpha\in\a{h}}.
$$
Then, the vertex operator map $Y(\cdot,x)$ is given by
\beqa
&&    Y(\alpha_1(-n_1) \cdots \alpha_j(-n_j){\bf 1},x) \no\\
    && \qquad = \: \frc1{(n_1-1)!} \lt( \frc{d}{dx} \rt)^{n_1-1}
\alpha_1(x) \cdots
    \frc1{(n_j-1)!} \lt( \frc{d}{dx} \rt)^{n_j-1} \alpha_j(x) \:
\eeqa
for $\alpha_k\in\a{h},\, n_k\in\Z_+,\; k=1,2,\ldots,j$, for all
$j\in\N$, where $\:\cdot\:$ is the usual normal ordering, which
brings $\alpha(n)$ with $n>0$ to the right. Choosing an
orthonormal basis $\{\b\alpha_q|q=1,\ldots,d\}$ of $\a{h}$, the
conformal vector is $\omega= \frc12 \sum_{q=1}^d \b\alpha_q(-1)
\b\alpha_q(-1){\bf 1}$. This implies in particular that the weight
of $\alpha(-n){\bf 1}$ is $n$:
\[
    L(0) \alpha(-n){\bf 1} = n\alpha(-n){\bf 1} \com{\alpha\in\a{h},\;
    n\in\Z_+}
\]
where we used
\[
    L(0) = \frc12 \sum_{n\in\frc1p\Z}\sum_{q=1}^d\,\:\b\alpha_q(n)\b\alpha_q(-n)\:\;.
\]
The isometry $\nu$ on $\a{h}$ lifts naturally to an automorphism
of the vertex operator algebra $S$, which we continue to call
$\nu$, of period $p$.

Then (cf. \cite{FLM2}), the various
properties of a vertex operator algebra are indeed satisfied by
the quadruplet $(V,Y,{\bf 1},\omega)$ just defined.

\ssect{Twisted modules}

We now proceed as in \cite{L1}, \cite{FLM1}, \cite{FLM2} and
\cite{DL2} to construct a space $S[\nu]$ that carries a natural
structure of $\nu$--twisted module for the vertex operator algebra
$S$. In these papers, the twisted module structure was constructed
assuming the minor hypothesis that $\nu$ preserves a rational lattice in $\a{h}$. We
show in Section \ref{sect-prooftmod} that the space $S[\nu]$ is a
twisted module, without the need for this minor assumption.

Consider a primitive $p$--th root of unity $\om$. For $r \in \Z$
set
$$
    \a{h}_{(r)} = \{\alpha\in\a{h} \;|\; \nu\alpha = \om^r\alpha\}
    \subset \a{h}.
$$
For $\alpha\in\a{h}$, denote by $\alpha_{(r)},\; r\in\Z$, its
projection on $\a{h}_{(r)}$. Define the $\nu$-twisted affine Lie
algebra $\hat{\a{h}}[\nu]$ associated with the abelian Lie algebra
$\a{h}$ by
\beq
    \hat{\a{h}}[\nu] = \coprod_{n\in\frc1p\Z} \a{h}_{(pn)} \otimes t^n
    \oplus \C C
\eeq
with
\beqa\label{commutrel}
    [\alpha\otimes t^m,\beta\otimes t^n] &=& \<\alpha,\beta\>
    m\delta_{m+n,0}\, C \com{\alpha\in\a{h}_{(pn)},\;
    \beta\in\a{h}_{(pm)},\; m,n\in\frc1p\Z} \no\\ {}
    [C,\hat{\a{h}}[\nu]] &=& 0 .
\eeqa
Set
\beq
    \hat{\a{h}}[\nu]^+ = \coprod_{n>0} \a{h}_{(pn)}\otimes t^n ,
    \qquad \hat{\a{h}}[\nu]^- = \coprod_{n<0} \a{h}_{(pn)} \otimes t^n.
\eeq
The subalgebra
\beq
    \hat{\a{h}}[\nu]^+ \oplus \hat{\a{h}}[\nu]^- \oplus \C C
\eeq
is a Heisenberg Lie algebra. Form the induced (level-one)
$\hat{\a{h}}[\nu]$-module
\beq
    S[\nu] = \mathcal{U}(\hat{\a{h}}[\nu])
    \otimes_{\mathcal{U}\lt( \hat{\a{h}}[\nu]^+
        \oplus \a{h}_{(0)} \oplus \C C\rt)} \C
        \simeq S(\hat{\a{h}}[\nu]^-) \for{(linearly),}
\eeq
where $\hat{\a{h}}[\nu]^+\oplus \a{h}_{(0)}$ acts trivially on
$\C$ and $C$ acts as 1; $\mathcal{U}(\cdot)$ denotes universal
enveloping algebra.  Then $S[\nu]$ is irreducible under the
Heisenberg algebra $\hat{\a{h}}[\nu]^+ \oplus \hat{\a{h}}[\nu]^-
\oplus \C C$. We will use the notation
$\alpha^\nu(n)\;(\alpha\in\a{h}_{(pn)},\, n\in\frc1p\Z)$ for the
action of $\alpha\otimes t^n \in \hat{\a{h}}[\nu]$ on $S[\nu]$.

\begin{remark}
The special case where $p=1$ ($\nu=1_{\a{h}}$) corresponds to the
$\hat{\a{h}}$-module $S$.
\end{remark}

The $\hat{\a{h}}[\nu]$-module $S[\nu]$ is naturally a
$\nu$--twisted module for the vertex operator algebra $S$. One
first constructs the following formal series acting on $S[\nu]$:
\beq\label{apar}
    \alpha^\nu(x) = \sum_{n\in\frc1p\Z} \alpha^\nu(n) x^{-n-1},
\eeq
as well as the formal series $W(v,x)$ for all $v\in S$:
\beqa
&&    W(\alpha_1(-n_1) \cdots \alpha_j(-n_j) {\bf 1},x) \no\\
&& \qquad    = \: \frc1{(n_1-1)!} \lt( \frc{d}{dx} \rt)^{n_1-1}
        \alpha_1^\nu(x) \cdots
        \frc1{(n_j-1)!} \lt( \frc{d}{dx} \rt)^{n_j-1} \alpha_j^\nu(x) \:
\eeqa
where $\alpha_k\in\a{h},\, n_k\in\Z_+,\; k=1,2,\ldots,j$, for all
$j\in\N$. The twisted vertex operator map $Y_{S[\nu]}(\cdot,x)$
acting on $S[\nu]$ is then given by
\beq\label{tvertex}
    Y_{S[\nu]}(v,x) = W(e^{\Delta_x} v,x) \com{v\in S}
\eeq where $\Delta_x$ is a certain formal operator involving the
formal variable $x$ \cite{FLM1}, \cite{FLM2}, \cite{DL2}. This
operator is trivial on $\alpha(-n){\bf 1}\in S\;(n\in\Z_+)$, so
that one has in particular \beq\label{tvertexalphax}
    Y_{S[\nu]}(\alpha(-n){\bf 1},x) = \frc1{(n-1)!}
        \lt( \frc{d}{dx} \rt)^{n-1}
    \alpha^\nu(x).
\eeq

One crucial (among others) role of the formal operator $\Delta_x$ is to make the
fixed--point subalgebra $\{u\,|\, \nu u=u\}$ act according to a
true module action. For instance, the conformal vector $\omega$ is
in the fixed--point subalgebra, so that the vertex operator
$Y_{S[\nu]}(\omega,x)$ generates a representation of the Virasoro
algebra on the space $S[\nu]$. This representation of the Virasoro
algebra was explicitly constructed in \cite{DL2}. As one can see
in the results of \cite{DL2} and as will become clear below, the
resulting representation $\Res_x\,xY_{S[\nu]}(\omega,x)$ of the
Virasoro generator $L(0)$ is not an (infinite) sum of
normal-ordered products the type
$\sum_{n\in\frc1p\Z}\,\:\alpha^\nu(n)\beta^\nu(-n)\:\;$; rather,
there is an extra term proportional to the identity on $S[\nu]$,
the so-called correction term, which appears because of the
operator $\Delta_x$. The correction term was calculated in
\cite{DL2} using the explicit action of $e^{\Delta_x}$ on
$\omega$. In the case of the period--2, $\nu=-1$ automorphism,
this action is given by \cite{FLM1}, \cite{FLM2}:
$$
    e^{\Delta_x}\omega = \omega + \frc1{16} (\dim h) x^{-2},
$$
and for general automorphism, the calculation was carried out in
\cite{DL2} (see also \cite{FFR} and \cite{FLM2}). These results
are relevant, for instance, in the construction of the moonshine
module \cite{FLM2}.

The calculation of the action of $\Delta_x$ on arbitrary elements
of $S$ is, however, a much more complicated task. Below we will
derive some identities among twisted vertex operators. One of the
important consequences of these identities, for us, will be to
give a tool to explicitly construct the twisted vertex operators
associated to elements of $S$ from the knowledge of the twisted
vertex operators associated to ``simpler'' elements, without
requiring the explicit knowledge of $\Delta_x$. In fact, these
identities allow us to construct recursively twisted vertex
operators associated to all elements of $S$ and to compute
$\Delta_x$, only starting from the knowledge that $\Delta_x$ is
trivial on $\alpha(-n){\bf 1}\in S\;(n\in\Z_+)$.

\sect{Commutativity and associativity properties}

This section follows closely similar sections of \cite{DLMi1} and
\cite{DLMi2}, and reproduces the results and some of the proofs.
We first recall the main commutativity and associativity
properties of vertex operators in the context of modules
(\cite{FLM2}, \cite{FHL}, \cite{DL1}, \cite{Li1}; cf. \cite{LL}),
and then we derive other identities somewhat analogous to these.
These other identities were stated and proven in \cite{DLMi2}, and
the most important ones were stated in \cite{DLMi1}. All these
identities will be generalized to twisted modules, still following
\cite{DLMi1,DLMi2}. Note that taking the module to be the vertex
operator algebra $V$ itself, the relations below specialize to
commutativity and associativity properties in vertex operator
algebras. We will give the proofs of the simplest identities only,
referring the reader to \cite{DLMi2} for all the proofs. Throughout this
and the next sections, we fix a vertex operator algebra $V$ and
a $V$-automorphism $\nu$ of period $p$, $\nu^p=1_V$.

\ssect{Formal commutativity and associativity for untwisted
modules}

We already stated the weak commutativity relation (\ref{wcom}) and
the weak associativity relation (\ref{wass}). They imply the main
``formal'' commutativity and associativity properties of vertex
operators, which, along with the fact that these properties are
equivalent to the Jacobi identity, can be formulated as follows
(see \cite{LL}):

\begin{theo} \label{theostructY}
Let $W$ be a vector space (not assumed to be graded) equipped with
a linear map $Y_W(\cdot,x)$ (\ref{vo}) such that the truncation
condition (\ref{truncation}) and the Jacobi identity
(\ref{jacobim}) hold. Then for $u,\,v\,\in V$ and $w\,\in W$,
there exist $k(u,v) \in\N$ and $l(u,w) \in \N$ and a (nonunique)
element $F(u,v,w;x_0,x_1,x_2)$ of $W((x_0,x_1,x_2))$ such that
\beqa
    &&  x_0^{k(u,v)} F(u,v,w;x_0,x_1,x_2) \in W[[x_0]]((x_1,x_2)), \no\\
    \label{propFuvw} &&  x_1^{l(u,w)} F(u,v,w;x_0,x_1,x_2) \in
    W[[x_1]]((x_0,x_2))
\eeqa
and
\beqa
    Y_W(u,x_1) Y_W(v,x_2)w &=& F(u,v,w;x_1-x_2,x_1,x_2), \no\\
    Y_W(v,x_2) Y_W(u,x_1)w &=& F(u,v,w;-x_2+x_1,x_1,x_2), \no\\
    \label{structY}
    Y_W(Y(u,x_0)v,x_2)w &=& F(u,v,w;x_0,x_2+x_0,x_2)
\eeqa
(where we are using the binomial expansion convention).
Conversely, let $W$ be a vector space equipped with a linear map
$Y_W(\cdot,x)$ (\ref{vo}) such that the truncation condition
(\ref{truncation}) and the statement above hold, except that
$k(u,v)$ ($\in\N$) and $l(u,w)$ ($\in \N$) may depend on all three
of $u,\,v$ and $w$. Then the Jacobi identity (\ref{jacobim})
holds.
\end{theo}

It is important to note that since $k(u,v)$ can be (and typically
is) greater than 0, the formal series $F(u,v,w;x_1-x_2,x_1,x_2)$
and $F(u,v,w;-x_2+x_1,x_1,x_2)$ are not in general equal. Along
with (\ref{propFuvw}), the first two equations of (\ref{structY})
represent formal commutativity, while the first and last equations
of (\ref{structY}) represent formal associativity, as formulated
in \cite{LL} (see also \cite{FLM2} and \cite{FHL}). The twisted
generalization of this theorem, written below, was proven in
\cite{DLMi2}.

\ssect{Additional relations in untwisted modules}

{}From the equations in Theorem \ref{theostructY}, we can derive a
number of relations similar to weak commutativity and weak
associativity but involving formal limit procedures (the meaning
of such formal limit procedures is recalled below). Even though
only one of these will be of use in the following sections, we
state here for completeness of the discussion the two relations
that are not ``easy'' consequences of weak commutativity and weak
associativity. These relations were proven in \cite{DLMi2}; we
report the proofs here.

The first relation can be expressed as follows:

\begin{theo}\label{theovwa}
With $W$ as in Theorem \ref{theostructY},
\begin{eqnarray} \label{varwass}
&&    \lim_{x_0 \to -x_2+x_1} \lt( (x_0+x_2)^{l(u,w)}
    Y_W(Y(u,x_0)v,x_2)w\rt) =
    x_1^{l(u,w)} Y_W(v,x_2) Y(u,x_1)w \nonumber \\
&&
\end{eqnarray}
for $u,\,v\,\in V$.
\end{theo}

The meaning of the formal limit
\begin{equation}\label{formallimit}
    \lim_{x_0 \to -x_2+x_1} \lt((x_0+x_2)^{l(u,w)}
    Y_W(Y(u,x_0)v,x_2)w\rt)
\end{equation}
is that one replaces each power of the formal variable $x_0$ in
the formal series $(x_0+x_2)^{l(u,w)} Y_W(Y(u,x_0)v,x_2)w$ by the
corresponding power of the formal series $-x_2+x_1$ (defined using
the binomial expansion convention). Notice again that the order of
$-x_2$ and $x_1$ is important in $-x_2+x_1$, according to the
binomial expansion convention.

{\em Proof of Theorem \ref{theovwa}:} Apply the limit $\lim_{x_0
\to -x_2+x_1}$ to the expression
\[
    (x_0+x_2)^{l(u,w)} Y_W(Y(u,x_0)v,x_2)w
\]
written as in the right--hand side of the third equation of
(\ref{structY}). This limit is well defined; indeed, the only
possible problems are the negative powers of $x_2+x_0$ in
$F(u,v,w,x_0,x_2+x_0,x_2)$, but they are cancelled out by the
factor $(x_0+x_2)^{l(u,w)}$. The resulting expression is read off
the second relation of (\ref{structY}). \eproof

\begin{remark}\label{rematrivialrelation}
It is instructive to consider the following relation, deceptively
similar to (\ref{varwass}), but that is in fact an immediate
consequence of weak associativity (\ref{wass}):
\beq \label{varwass2}
    \lim_{x_0 \to x_1-x_2} \lt( (x_0+x_2)^{l(u,w)}
    Y_W(Y(u,x_0)v,x_2)w\rt) =
    x_1^{l(u,w)} Y_W(u,x_1) Y_W(v,x_2)w.
\eeq
More precisely, it can be obtained by noticing that the
replacement of $x_0$ by $x_1-x_2$ independently in each factor in
the expression as written on the left--hand side of (\ref{wass})
is well defined. We emphasize that, by contrast, the relation
(\ref{varwass}) {\em cannot} be obtained in such a manner. Indeed,
although the formal limit procedure $\lim_{x_0 \to -x_2+x_1}$ is
of course well defined on the series on both sides of
(\ref{wass}), one cannot replace $x_0$ by $-x_2+x_1$ either in the
factor $Y_W(u,x_0+x_2)Y_W(v,x_2)w$ on the left--hand side or in
the factor $Y_W(Y(u,x_0)v,x_2)w$ on the right--hand side of
(\ref{wass}).
\end{remark}

The second nontrivial relation, which we call {\em modified weak
associativity}, will be important when generalized to twisted
modules. It was first written in \cite{DLMi1}. It is stated as:

\begin{theo}\label{theomwa}
With $W$ as in Theorem \ref{theostructY},
\beq\label{modwass}
    \lim_{x_1 \to x_2+x_0} \lt( (x_1-x_2)^{k(u,v)} Y_W(u,x_1)
    Y_W(v,x_2) \rt) =
    x_0^{k(u,v)} Y_W(Y(u,x_0)v,x_2)
\eeq
for $u,\,v\,\in V$.
\end{theo}

\proof Apply the limit $\lim_{x_1 \to x_2+x_0}$ to the expression
\[
    (x_1-x_2)^{k(u,v)} Y_W(u,x_1) Y_W(v,x_2)
\]
written as in the right--hand side of the first equation of
(\ref{structY}). This limit is well defined, since negative powers
of $x_1-x_2$ in $F(u,v,w;x_1-x_2,x_1,x_2)$ are cancelled out by
the factor $(x_1-x_2)^{k(u,v)}$. The resulting expression is read
off the third relation of (\ref{structY}). \eproof

\begin{remark}\label{remares}
Equation (\ref{modwass}) can be written in the following form:
\beq
    \lim_{x_1 \to x_2+x_0} \lt( \lt(\frc{x_1-x_2}{x_0}\rt)^{k(u,v)}
    Y_W(u,x_1) Y_W(v,x_2) \rt) =
    Y_W(Y(u,x_0)v,x_2).
\eeq
The factor $\lt(\frc{x_1-x_2}{x_0}\rt)^{k(u,v)}$ appearing in
front of the product of two vertex operators on the left--hand
side is crucial in giving a well--defined limit, but when the
limit is applied to this factor without the product of vertex
operators, the result is simply $1$.  We will call such a factor a
``resolving factor''. Its power is apparent, in particular, in the
proof of the main commutator formula (6.1) of \cite{DLMi2}: it
allows one to evaluate nontrivial limits of sums of terms with
cancelling ``singularities'' in a straightforward fashion,
evaluating the limit of each term independently. Its power will
also be clear, in the present paper, when constructing the twisted
vertex operator map $Y_{S[\nu]}(\cdot,x)$ and when studying the
algebra $\h{D}^+$ defined in Section \ref{sect-D}.
\end{remark}

\ssect{Formal commutativity and associativity for twisted modules}

We derive below various commutativity and associativity properties
of twisted vertex operators. In order to express some of these
properties, we need one more element of formal calculus: a certain
projection operator (see \cite{DLMi2}). Consider the operator
$P_{[[x_0,x_0^{-1}]]}$ acting on the space $\mathbb{C}\{x_0\}$ of
formal series with any complex powers of $x_0$, which projects to
the formal series with integral powers of $x_0$:
\beq\label{projection}
    P_{[[x_0,x_0^{-1}]]} \;:\;
    \mathbb{C}\{x_0\} \to \mathbb{C}[[x_0,x_0^{-1}]] \;.
\eeq
We will extend the meaning of this notation in the obvious way to
projections acting on formal series with coefficients lying in
vector spaces other than $\C$, vector spaces which might
themselves be spaces of formal series in other formal variables.
Notice that when this projection operator acts on a formal series
in $x_0$ with powers that are in $\frc1p \Z$, for instance on
$f(x_0) \in \C [[x_0^{1/p},x_0^{-1/p}]]$, it can be described by
an explicit formula:
\[
    P_{[[x_0,x_0^{-1}]]} f(x_0) = \frc1p \sum_{r=0}^{p-1}
    \lt( \lim_{x^{1/p}\to \om^r x_0^{1/p}} f(x) \rt).
\]
(See Remark \ref{remaformallimitfrac} below for the meaning of
formal limit procedures involving fractional powers of formal
variables.)  We will also extend this projection notation to
different kinds of formal series in obvious ways.  For instance,
\[
    P_{x_0^{q/p} [[x_0,x_0^{-1}]]} \;:\;
    \mathbb{C}\{x_0\} \to \mathbb{C}\,x_0^{q/p}[[x_0,x_0^{-1}]] \;.
\]
Again, of course, we will extend the meaning of this notation
to formal series with coefficients in vector spaces other than
$\C$.

The twisted Jacobi identity (\ref{jacobitm}) implies twisted
versions of weak commutativity and weak associativity ($u,\,v\,\in
V,\;w\in M$):
\beqa \label{wcomt}
    &&      (x_2-x_1)^{k} Y_M(v,x_2)Y_M(u,x_1) =
    (x_2-x_1)^{k} Y_M(u,x_1) Y_M(v,x_2)
    \nn && \\
    \label{wasst}
    && P_{[[x_0,x_0^{-1}]]} \lt( (x_0+x_2)^{l}
    Y_M(u,x_0+x_2)Y_M(v,x_2)w \rt) \no\\
    &&      \qquad = (x_2+x_0)^{l} \frc1p
    \sum_{r=0}^{p-1}
    \omega^{-lrp}_p Y_M(Y(\nu^ru,x_0)v,x_2)w.
\eeqa
These relations are valid for all large enough $k \in \N$ and $l
\in \frc1p \N$, their minimum value depending respectively on
$u,\,v$ and on $u,\,w$. For definiteness, we will denote these
minimum values by $k(u,v)$ and $l(u,w)$, respectively (they depend
also on the module $M$; in particular, they differ from the
integer numbers $k(u,v)$ and $l(u,w)$ used in the previous
subsection in connection with the module $W$). As in the untwisted
case, these relations imply the main ``formal'' commutativity and
associativity properties of twisted vertex operators \cite{Li2},
which, along with with the fact that these properties are
equivalent to the Jacobi identity, can be formulated as follows
(it was first formulated in this form in \cite{DLMi1}):

\begin{theo} \label{theostructYt}
Let $M$ be a vector space (not assumed to be graded) equipped with
a linear map $Y_M(\cdot,x)$ (\ref{tvo}) such that the truncation
condition (\ref{truncationt}) and the Jacobi identity
(\ref{jacobitm}) hold. Then for $u,\,v\,\in V$ and $w\,\in M$,
there exist $k(u,v) \in\N$ and $l(u,w) \in \frc1p\N$ and a
(nonunique) element $F(u,v,w;x_0,x_1,x_2)$ of
$M((x_0,x_1^{1/p},x_2^{1/p}))$ such that
\beqa
    &&  x_0^{k(u,v)} F(u,v,w;x_0,x_1,x_2)
    \in M[[x_0]]((x_1^{1/p},x_2^{1/p})), \no\\
    \label{propFuvwt}
    &&  x_1^{l(u,w)} F(u,v,w;x_0,x_1,x_2) \in
    M[[x_1^{1/p}]]((x_0,x_2^{1/p}))
\eeqa
and
\beqa
    Y_M(u,x_1) Y_M(v,x_2)w &=& F(u,v,w;x_1-x_2,x_1,x_2), \no\\
    Y_M(v,x_2) Y_M(u,x_1)w &=& F(u,v,w;-x_2+x_1,x_1,x_2), \no\\
    \label{structYt}
    Y_M(Y(\nu^{-s}u,x_0)v,x_2)w &=&
    \lim_{x_1^{1/p} \to \om^s (x_2+x_0)^{1/p}}
     F(u,v,w;x_0,x_1,x_2)
\eeqa
for $s\in\Z$ (where we are using the binomial expansion
convention). Conversely, let $M$ be a vector space equipped with a
linear map $Y_M(\cdot,x)$ (\ref{tvo}) such that the truncation
condition (\ref{truncation}) and the statement above hold, except
that $k(u,v)$ ($\in\N$) and $l(u,w)$ ($\in \frc1p\N$) may depend
on all three of $u,\,v$ and $w$. Then the Jacobi identity
(\ref{jacobitm}) holds.
\end{theo}

This theorem, as well as (\ref{wcomt}) and (\ref{wasst}), were
proven in \cite{DLMi2}.

\begin{remark}\label{remaformallimitfrac}
Formal limit procedures involving fractional powers of formal
variables like $x_1^{1/p}$ have the same meaning as in
(\ref{formallimit}), but with $x_1^{1/p}$ being treated as a
formal variable by itself. For instance, the formal limit
procedure
\[
    \displaystyle{\lim_{x_1^{1/p} \to \om^s
    (x_2+x_0)^{1/p}}F(u,v,w;x_0,x_1,x_2)}
\]
above means that one replaces each integral power of the formal
variable $x_1^{1/p}$ in the formal series $F(u,v,w;x_0,x_1,x_2)$
by the corresponding power of the formal series $\om^s
(x_2+x_0)^{1/p}$ (defined using the binomial expansion
convention).
\end{remark}

\begin{remark}
Note that this theorem, and in particular its proof in
\cite{DLMi2}, illustrates the phenomenon, which arises again and
again throughout the theory of vertex operator algebras, that
formal calculus inherently involves just as much ``analysis'' as
``algebra'': in many relations there are integers that can be left
unspecified, except for their minimum values, and the proof
involves taking these integers ``large enough''. Recall that
essentially the same issues arose for example in the use of formal
calculus for the proof of the Jacobi identity for (twisted) vertex
operators in \cite{FLM2} (see Chapters 8 and 9). This is certainly
not surprising, since we are using the Jacobi identity (for all
twisting automorphisms) in order to prove, in a different
approach, properties of (twisted) vertex operators.
\end{remark}

Along with (\ref{propFuvwt}), the first two equations of
(\ref{structYt}) represent what we call {\em formal commutativity}
for twisted vertex operators, while the first and last equations
of (\ref{structYt}) represent {\em formal associativity} for
twisted vertex operators. When specialized to the untwisted case
$p=1$ ($\nu=1_V$), these two relations lead respectively to the
usual formal commutativity and formal associativity for vertex
operators, as described in (\ref{structY}).

\ssect{Additional relations in twisted modules}

As in the case of ordinary vertex operators, one can write other
relations involving formal limit procedures. These relations were
proven in \cite{DLMi2}; we report the proofs here. Among them, two
cannot be directly obtained {}from weak commutativity and weak
associativity. One of these, the relation generalizing
(\ref{varwass}), is stated as follows:

\begin{theo}\label{theovwat}
With $M$ as in Theorem \ref{theostructYt},
\beqa &&
    \lim_{x_0 \to -x_2+x_1} \lt( (x_2+x_0)^{l} \frc1p
    \sum_{r=0}^{p-1}
    \omega^{-lrp}_p Y_M(Y(\nu^ru,x_0)v,x_2)w \rt) \no\\
    \label{varwasst}
    &&  \qquad =  P_{[[x_1,x_1^{-1}]]} \lt(
        x_1^{l} Y_M(u,x_2)Y_M(v,x_1)w \rt),
\eeqa
for all $l\in\frc1p\Z,\; l\ge l(u,w)$.
\end{theo}

\proof This is proved along the lines of the proof of Theorem
\ref{theovwa}, with some additions due to the fractional powers.
One uses the third equation of (\ref{structYt}) in order to
rewrite the left--hand side of (\ref{varwasst}) as
\[
    \lim_{x_0 \to -x_2+x_1} \lt( (x_2+x_0)^{l} \frc1p \sum_{r=0}^{p-1}
    \omega^{-lrp}_p \lim_{x_3^{1/p} \to
    \om^{-r}(x_2+x_0)^{1/p}}
    F(u,v,w;x_0,x_3,x_2) \rt).
\]
The sum over $r$ keeps only the terms in which $x_2+x_0$ is raised
to a power which has a fractional part equal to the negative of
the fractional part of $l$. Multiplying by $(x_2+x_0)^l$, for any
$l\in\frc1p\Z,\; l\ge l(u,w)$, brings the remaining series to a
series with finitely many negative powers of $x_2$ (as well as
$x_0$), to which it is possible to apply the limit $\lim_{x_0 \to
-x_2+x_1}$. This limit of course brings only integer powers of
$x_1$, and the right--hand side of (\ref{varwasst}) can be
obtained {}from the second equation of (\ref{structYt}). \eproof

\begin{remark}
A relation similar to the last one, but that is a direct
consequence of weak associativity (\ref{wasst}), is
\beqa
    && \lim_{x_0 \to x_1-x_2} \lt( (x_2+x_0)^{l} \frc1p
        \sum_{r=0}^{p-1} \omega^{-lrp}_p Y_M(Y(\nu^ru,x_0)v,x_2)w
        \rt) \no\\
    \label{varwass2t}
    && \qquad= P_{[[x_1,x_1^{-1}]]} \lt( x_1^{l}
    Y_M(u,x_1)Y_M(v,x_2)w \rt)
\eeqa
for all $l\in\frc1p\Z,\; l\ge l(u,w)$. This generalizes
(\ref{varwass2}) (see the comments in Remark
\ref{rematrivialrelation}). It can be obtained by applying the
formal limit involved in the left--hand side to both sides of
(\ref{wasst}).
\end{remark}

The most important relation for our purposes, which was first
stated in \cite{DLMi1}, generalizing (\ref{modwass}) and which we
call {\em modified weak associativity} for twisted vertex
operators, is given by the following theorem:

\begin{theo}\label{theomwat}
With $M$ as in Theorem \ref{theostructYt},
\beqa
    &&  \lim_{x_1^{1/p} \to \om^{s}(x_2+x_0)^{1/p}}
    \lt( (x_1-x_2)^{k(u,v)} Y_M(u,x_1) Y_M(v,x_2) \rt) \no\\
    \label{modwasst} && \qquad\qquad= x_0^{k(u,v)}Y_M(Y(\nu^{-s}u,x_0)v,x_2)
\eeqa
for $u,\,v\,\in V$ and $s\in \Z$.
\end{theo}

\proof The proof is a straightforward generalization of the proof
of Theorem \ref{theomwa}. \eproof

\begin{remark}
The specialization of Theorems \ref{theostructYt} and
\ref{theomwat} to the untwisted case $p=1$ and $M=W$ gives,
respectively, Theorems \ref{theostructY} and \ref{theomwa}.
\end{remark}

Finally, we derive a simple relation, proved in \cite{Li2}, that specifies the structure
of the formal series $Y_M(u,x)$.

\begin{theo}\label{theotransnu}
With $M$ as in Theorem \ref{theostructYt},
\beq \label{transnu}
    \lim_{x_1^{1/p} \to \om^s x^{1/p}} Y_M(\nu^{s}u,x_1) = Y_M(u,x)
\eeq
for $u\in V$ and $s\in \Z$.
\end{theo}

\proof In the Jacobi identity (\ref{jacobitm}), replace $u$ by
$\nu^s u$ and $x_1^{1/p}$ by $\om^s x^{1/p}$. The right--hand side
becomes
\[
    \frc1p x_2^{-1} \sum_{r=0}^{p-1}
    \delta\lt(\om^{r+s}\lt(\frc{x-x_0}{x_2}\rt)^{1/p}\rt)
        Y_M(Y(\nu^{r+s} u,x_0)v,x_2),
\]
which is independent of $s$, as is apparent if we make the shift
in the summation variable $r\mapsto r-s$. Hence the left--hand
side is also independent of $s$. Choosing $v={\bf 1}$ and using
the vacuum property (\ref{vacuumt}), this gives
\beqa
    &&    \lt( x_0^{-1}\dpf{x-x_2}{x_0}
    - x_0^{-1}\dpf{x_2-x}{-x_0}
    \rt) \lim_{x_1^{1/p} \to \om^s x^{1/p}} Y_M(\nu^s
    u,x_1) \\
    &&\qquad\qquad\qquad    =
    \lt( x_0^{-1}\dpf{x-x_2}{x_0}
    - x_0^{-1}\dpf{x_2-x}{-x_0}
    \rt) Y_M(u,x) \no
\eeqa
which, upon using (\ref{3delta}) and taking $\Res_{x_2}$, gives
(\ref{transnu}). \eproof

{}From Theorem \ref{theotransnu}, we directly have the following
corollary:

\begin{corol}
With $M$ as in Theorem \ref{theostructYt},
\[
    Y_M(u,x) = \sum_{n\in\Z+q/p} u_nx^{-n-1} \for{for} u\in V,\ \nu u = \om^q
    u,\ q\in\Z.
\]
\end{corol}

\sect{Equivalence and construction theorems from modified weak
associativity}

This section presents new results related to modified weak
associativity. We will show, loosely speaking, that modified weak
associativity (\ref{modwasst}) and weak commutativity
(\ref{wcomt}) are equivalent to the Jacobi identity
(\ref{jacobitm}) (``equivalence theorem''). Then we will show that if the modified weak
associativity for twisted vertex operators is valid for all pairs
$u,w$ (to be put in (\ref{modwasst}) instead of the ordered pair
$u,v$) with $u\in U$ and $w\in V$, where $U$ is a
generating subset of a vertex operator
algebra $V$, and if weak commutativity for twisted vertex
operators is valid for all pairs $u,v$ with $u,v\,\in U$, then
both modified weak associativity and weak commutativity hold for
the whole vertex operator algebra $V$ (``construction theorem'').  Similar construction theorems
were proved by Li in the untwisted and twisted cases \cite{Li1,Li2} (cf. \cite{LL}), using the powerful idea of ``local systems
of (twisted) vertex operators.'' Here we start from similar ideas but we make use of modified weak associativity discovered in \cite{DLMi1,DLMi2}, in order to illustrate one of its applications. We find it instructive to give a direct proof of our construction theorem, although there may be a shorter route from Li's construction theorems. In the next section, our two theorems will allow us to construct in a relatively simple way the twisted vertex operator map for $S[\nu]$ -- in particular, we will show how the use of modified weak associativity gives a new explicit form for the operator $\Delta_x$ -- and to
prove the twisted module structure for $S[\nu]$. We recall that throughout this section, $V$ is a vertex operator algebra and $\nu$ is an automorphism of $V$.

\ssect{Equivalence theorem}

It is well known in the theory of vertex operator algebras that,
under natural conditions, the weak commutativity relation
(\ref{wcom}) and the weak associativity relation (\ref{wass}) for
untwisted modules are equivalent to the Jacobi identity
(\ref{jacobim}). It is a simple matter to show that this statement
is also true when weak associativity is replaced by modified weak
associativity. We state this more generally for twisted modules
(and for the twisted Jacobi identity (\ref{jacobitm})) in the
following theorem.

\begin{theo} \label{theoequil}
Let $M$ be a vector space (not assumed to be graded) equipped with
a linear map $Y_M(\cdot,x)$ (\ref{tvo}) such that the truncation
condition (\ref{truncationt}) hold.

If modified weak associativity (\ref{modwasst}) holds, for some
$k'(u,v)\in\N$, and weak commutativity (\ref{wcomt}) holds, for a
possibly different $k(u,v)\in\N$, then we may change $k'(u,v)$ (in particular, we may lower it)
to $k'(u,v)=k(u,v)$, and the twisted Jacobi identity (\ref{jacobitm})
holds. Also, the vector space $M$ is a twisted module if,
additionally, $M$ is ${\mathbb Q}-graded$ and quasi-finite as in the first lines of
Definition \ref{VOAmodulet}, with the $L(0)$-weight property (\ref{L0weightmt}),
and the vacuum property (\ref{vacuumt}) holds.

On the other hand, if the twisted Jacobi identity (\ref{jacobitm})
holds, then modified weak associativity (\ref{modwasst}) and weak
commutativity (\ref{wcomt}) hold.
\end{theo}

\begin{rema} Note that this theorem can be specialized to $p=1$,
applying then to modules $W$. It can also be used to show that
some vector space $V$ is a vertex operator algebra; more precisely,
in the definition \ref{VOA}, the Jacobi identity can be replaced by modified
weak associativity and weak commutativity.
\end{rema}

{\em Proof of Theorem \ref{theoequil}.} The last sentence of the
theorem was proven already by proving modified weak
associativity and weak commutativity above.

Weak commutativity (\ref{wcomt}), when both sides are applied on
an element $w$ of $M$, and the truncation condition immediately
imply the first equation of (\ref{propFuvwt}), and the first two
equations of (\ref{structYt}). Then, the limit on the
left-hand side of (\ref{modwasst}) with $k(u,v)$ replaced by
$k''(u,v)$ certainly exists for all $k''(u,v)\ge
k(u,v)$, and the result is the same for any $k''(u,v)$ that makes the limit exist, up to the obvious power of $x_0$
(because if all limits exist, then the product of the limits is the limit of the product).
Since we know that this limit gives the right-hand side for
some $k'(u,v)$, we may change it (and lower it if necessary) to $k'(u,v)=k(u,v)$. Then, writing the
product of twisted vertex operators on the left-hand side of
(\ref{modwasst}) as on the right-hand side of the first equation
of (\ref{structYt}), we see that (\ref{modwasst}) implies the
third equation of (\ref{structYt}). Hence, by Theorem
\ref{theostructYt}, the twisted Jacobi identity holds. With the
additional conditions stated in the theorem, all other parts of
the definition (\ref{VOAmodulet}) are satisfied (in particular,
the fact that $V$ is a vertex operator algebra implies that the
Virasoro commutation relations are also satisfied, and the $L(-1)$-derivative
property for vertex operator algebra (\ref{Lm1der}) implies, using (\ref{modwasst}) with $v={\bf 1}$,
the corresponding property for twisted modules (\ref{Lm1dermt})) and $M$ is a twisted
module for $V$. \eproof

\ssect{Construction theorem}

Consider a generating subset $U\subset V$ of the vertex operator
algebra $V$, defined as follows.
\begin{defi} \label{defgen}
A generating subset $U\subset V$ of a vertex operator algebra is
a subset such that all elements of $V$ can be written as linear combinations of
elements of the form $u_mv_n\cdots {\bf 1}$ for $u,v,\ldots \in U$ and
$m,n,\ldots\in\Z$.
\end{defi}
Suppose we are able to define a twisted
vertex operator map $Y_M(\cdot,x)$ as in (\ref{tvo}) such that the truncation condition (\ref{truncationt}) holds for all $v\in V$,
weak commutativity (\ref{wcomt}) holds for all $u,v\,\in U$, and modified weak associativity (\ref{modwasst}) holds for all $u\in U$ and $v\in V$.
Theorem \ref{theoconstr1} below along with Theorem \ref{theoequil} tell us that this is enough to have the twisted Jacobi identity, and, with additional mild
conditions (see Theorem \ref{theoequil}) to have a twisted module.

\begin{rema} Theorem \ref{theoconstr1} of course requires us to already have a twisted vertex operator map $Y_M(\cdot,x)$ on the full vertex operator algebra $V$.
On the other hand, once we have the map on $U$ only (with the properties above),
it is easy to extend it to a map on the set of symbols of the type $u_m v_n \cdots w$ for $u,v,\ldots,w \in U$ and $m,n,\ldots\in\Z$
by recursive use of modified weak associativity (\ref{modwasst}) with $u\in U$ and $v$ an element in this set of symbols. However, the vertex operator
algebra $V$ is the span of this set of symbols {\em with linear relations amongst them}; these relations depend on the particular vertex operator algebra
at hand. In order to have a twisted module, it is essential to verify that the
map on this set of symbols is well-defined on $V$; that is, that it is in agreement with these linear relations. This is extremely nontrivial, and there does not
seem to be yet a general theorem as to when that happens.
Whatever these relations are, they have to imply those coming from the Jacobi identity (since $V$ is a vertex operator algebra). 
It is possible to show that all linear relations coming from the Jacobi identity are satisfied by this construction. This is beyond the scope of the present
paper, but we hope to clarify some of these issues in a future work.
\end{rema}

Note again that the theorem below and the issues discussed in the remark above are in close relation with results of Li \cite{Li2}; cf. \cite{LL}.

\begin{theo} \label{theoconstr1}
Let $M$ be a vector space (not assumed to be graded) equipped with
a linear map $Y_M(\cdot,x)$ (\ref{tvo}) such that the truncation
condition (\ref{truncationt}) holds. Fix a generating subset $U\subset
V$, as defined in Definition \ref{defgen}.

If weak commutativity (\ref{wcomt}) is satisfied for all $u,v\,\in
U$ (and fix $k(u,v)\in\N$ to be the lowest integer that can be
taken both in (\ref{wcomt}) and in weak commutativity for the
vertex operator algebra $V$), and if modified weak associativity
(\ref{modwasst}) is satisfied for all $u\in U$ and for all $v\in
V$ (for some $k'(u,v)\in\N$ that may be different from $k(u,v)$),
then both weak commutativity and modified weak associativity are
satisfied for all $u,v\,\in V$. Further, we may change $k'(u,v)$
(in particular, we may lower it) to $k'(u,v)=k(u,v)$ for $u,v\,\in V$, and these integers may be
taken to satisfy the formula
\beq\label{kdesc}
    k(u_{n_0+m}v,w) = k(w,u) + k(w,v) + m
\eeq
for all $u,v,w\,\in V$ where $n_0$ is the highest integer such
that $u_{n_0}v\neq 0$. The same integers may also be taken in weak
commutativity for the vertex operator algebra $V$.
\end{theo}

\proof We start by showing weak commutativity (\ref{wcomt}), and
the equation (\ref{kdesc}). Under the assumptions of the theorem,
we have, for $u,v,w\,\in U$,
\beqa
    &&  (x-x_2-x_0)^{k(w,u)} (x-x_2)^{k(w,v)} Y_M(w,x)\cdot
    \no\\ && \qquad\qquad \cdot
    \lim_{x_1^{1/p} \to \om^{s}(x_2+x_0)^{1/p}}
    \lt( (x_1-x_2)^{k(u,v)} Y_M(u,x_1) Y_M(v,x_2) \rt)
    \no\\ && =
    \lim_{x_1^{1/p} \to \om^{s}(x_2+x_0)^{1/p}}
    \lt(  (x-x_1)^{k(w,u)} (x-x_2)^{k(w,v)} Y_M(w,x) \cdot
    \rt. \no\\ && \qquad\qquad \cdot \lt.
    (x_1-x_2)^{k(u,v)} Y_M(u,x_1) Y_M(v,x_2) \rt)
    \no\\ && =
    \lim_{x_1^{1/p} \to \om^{s}(x_2+x_0)^{1/p}}
    \lt( (x_1-x_2)^{k(u,v)} Y_M(u,x_1) Y_M(v,x_2) \cdot
    \rt. \no\\ && \qquad\qquad \cdot \lt.
    (x-x_1)^{k(w,u)} (x-x_2)^{k(w,v)} Y_M(w,x) \rt)
    \no\\  && =
    \lim_{x_1^{1/p} \to \om^{s}(x_2+x_0)^{1/p}}
    \lt( (x_1-x_2)^{k(u,v)} Y_M(u,x_1) Y_M(v,x_2)  \rt)
    \cdot \no\\ && \qquad\qquad \cdot
    (x-x_2-x_0)^{k(w,u)} (x-x_2)^{k(w,v)} Y_M(w,x)
\eeqa
hence
\beqa\label{prwc}
&&    (x-x_2-x_0)^{k(w,u)} (x-x_2)^{k(w,v)} \,
    Y_M(w,x) Y_M(Y(\nu^{-s}u,x_0)v,x_2)
    \\ && \qquad\qquad =
    (x-x_2-x_0)^{k(w,u)} (x-x_2)^{k(w,v)}
    Y_M(Y(\nu^{-s}u,x_0)v,x_2) Y_M(w,x)~. \no
\eeqa
Both sides have finitely many negative powers of $x_0$. Take the
lowest power:
\beqa\label{prwc2}
&&    (x-x_2)^{k(w,u)+k(w,v)} \,
    Y_M(w,x) Y_M((\nu^{-s}u)_{n_0}v,x_2)
    \\ && \qquad\qquad =
    (x-x_2)^{k(w,u)+k(w,v)}
    Y_M((\nu^{-s}u)_{n_0}v,x_2) Y_M(w,x)~. \no
\eeqa
This proves weak commutativity for all pairs $u_{n_0}v,w$ with
$u,v,w\,\in U$, and we may take $k(u_{n_0}v,w) = k(w,u) + k(w,v)$
(although we have not proven that this value is the minimum one
for which weak commutativity is valid), where $n_0$ is the highest
integer such that $u_{n_0}v\neq 0$.

The next power of $x_0$ in the equation (\ref{prwc}) contains, on
each side, two terms of the same type as those above: one with
$(\nu^{-s}u)_{n_0+1}v$, the other with $(\nu^{-s}u)_{n_0} v$. The
terms with $(\nu^{-s}u)_{n_0}v$ are multiplied by
$(x-x_2)^{k(w,u)+k(w,v)-1}$, and those with
$(\nu^{-s}u)_{n_0+1}v$, by $(x-x_2)^{k(w,u)+k(w,v)}$. Multiplying
the resulting equation through by $x-x_2$, we can use the result
(\ref{prwc2}) and we obtain weak commutativity for all pairs
$u_{n_0+1}v,w$ with $u,v,w\,\in U$, where $n_0$ is again highest
integer such that $u_{n_0}v\neq 0$. We can take $k(u_{n_0+1}v,w) =
k(w,u) + k(w,v)+1$. Repeating the process, we obtain weak
commutativity for all pairs $u_nv,w$ with $u,v,w\,\in U$ and for
all $n\in \Z$, with $k(u_{n_0+m}v,w) = k(w,u) + k(w,v)+m$.
Replacing $v$ by $u_nv$ for $u,v\,\in U$ in the argument above,
and repeating, we obtain weak commutativity for all pairs $v,w$
with $v\in V$ and $w\in U$. Finally, we can repeat the full
argument above with $u,v\in U$ and $w\in V$. The induction from
that gives us weak commutativity for all pairs $v,w$ with $v\in V$
and $w\in V$, and in particular gives us (\ref{kdesc}).

Note that the same arguments can be used for the vertex operator
$Y(\cdot,x)$, so that equation (\ref{kdesc}) may be taken to hold
also for the $k(u,v)$ involved in weak commutativity for the
vertex operator algebra $V$. This implies that the integers
$k(u,v)$ involved in weak commutativity for $Y_M(\cdot,x)$ and
those involved in weak commutativity for the vertex operator
algebra $V$ may both be taken to be given by (\ref{kdesc}) for all elements of the
vertex operator algebra. We will do that for the rest of the
proof.

Note also that weak commutativity (\ref{wcomt}) implies, through the arguments
of the proof of Theorem \ref{theoequil}, that in
modified weak associativity (\ref{modwasst}) for $u\in U$ and
$v\in V$, which holds by assumption with some integer $k'(u,v)$,
we can change $k'(u,v)$ to $k'(u,v)=k(u,v)$. We will also do that
for the rest of the proof.

Next we show modified weak associativity. First note the following
lemma.

\begin{lemma} \label{lemma3oper}
For $u,v,w\,\in V$ and $\phi\in M$, there exists a (nonunique)
element $F(u,v,w,\phi;x_0,x_4,x_5,x_1,x_2,x_3)$ of
$M((x_0,x_4,x_5,x_1^{1/p},x_2^{1/p},x_3^{1/p}))$ such that
\beqa
    x_0^{k(u,v)}F(u,v,w,\phi;x_0,x_4,x_5,x_1,x_2,x_3) \in
    M[[x_0]]((x_4,x_5,x_1^{1/p},x_2^{1/p},x_3^{1/p})),\no\\
    x_4^{k(u,w)}F(u,v,w,\phi;x_0,x_4,x_5,x_1,x_2,x_3) \in
    M[[x_4]]((x_0,x_5,x_1^{1/p},x_2^{1/p},x_3^{1/p})),\no\\
    x_5^{k(v,w)}F(u,v,w,\phi;x_0,x_4,x_5,x_1,x_2,x_3) \in
    M[[x_5]]((x_0,x_4,x_1^{1/p},x_2^{1/p},x_3^{1/p}))
\eeqa
and
\beq
    Y_M(u,x_1)Y_M(v,x_2)Y_M(w,x_3)\phi =
    F(u,v,w,\phi;x_1-x_2,x_1-x_3,x_2-x_3,x_1,x_2,x_3)
\eeq
where $k(u,v),k(u,w),k(v,w)\in\N$ can be taken as those that
appear in the weak commutativity relation (\ref{wcomt}).
\end{lemma}

This is an immediate consequence of weak commutativity
(\ref{wcomt}). Indeed, consider
\[
    (x_1-x_2)^{k(u,v)} (x_1-x_3)^{k(u,w)} (x_2-x_3)^{k(v,w)}
    Y_M(u,x_1)Y_M(v,x_2)Y_M(w,x_3)\phi~.
\]
Thanks to the factor $(x_1-x_2)^{k(u,v)} (x_1-x_3)^{k(u,w)}
(x_2-x_3)^{k(v,w)}$, the twisted vertex operators may be written
in any order. Looking at different orders and using the truncation
property, we see the expression has finitely many negative powers
of $x_3$, finitely many negative powers of $x_2$ and finitely many
negative powers of $x_1$. This shows the lemma.

Now, take $u,v,\in U$ and $w\in V$, and again $\phi\in M$. From
modified weak associativity for vertex operator algebra
(\ref{modwass}) (recall that the module $W$ in that equation may
be replaced by the vertex operator algebra itself $V$), which we
will need only for the pair $u,v$, and from modified weak
associativity for twisted modules, which is valid by assumption
when, in (\ref{modwasst}), the (ordered) pair $u,v$ is replaced by
the pair $v,w$ as well as when it is replaced by the pair
$u,v_nw$, we have
\beqa
    &&    Y_M(Y(Y(u,x_0)v,x_5)w,x_3)\phi \no \\ && =
    \lim_{x_4\to x_5+x_0}
    \lt(\frc{x_4-x_5}{x_0}\rt)^{k(u,v)} Y_M(Y(u,x_4)Y(v,x_5)w,x_3)
    \no \\ &&=
    \lim_{x_4\to x_5+x_0}
    \lt(\frc{x_4-x_5}{x_0}\rt)^{k(u,v)}
    \sum_{n\in\Z} x_5^{-n-1} Y_M(Y(u,x_4)v_nw,x_3)\phi
    \no\\ &&=
    \lim_{x_4\to x_5+x_0}
    \lt(\frc{x_4-x_5}{x_0}\rt)^{k(u,v)}
    \sum_{n\in\Z} x_5^{-n-1}
    \lim_{x_1^{\frc1p}\to (x_2+x_4)^{\frc1p}}
    \lt(\frc{x_1-x_3}{x_4}\rt)^{k(u,v_nw)} \cdot \no\\ &&\qquad
    \cdot Y_M(u,x_1)Y_M(v_nw,x_3)\phi \no
    \no\\ &&=
    \lim_{x_4\to x_5+x_0}
    \lt(\frc{x_4-x_5}{x_0}\rt)^{k(u,v)}
    \sum_{n\in\Z}\quad
    \lim_{x_1^{\frc1p}\to (x_3+x_4)^{\frc1p}}
    \lt(\frc{x_1-x_3}{x_4}\rt)^{k(u,v_nw)} \cdot \no\\ &&\qquad
    \cdot P_{x_5^{-n-1}} Y_M(u,x_1)Y_M(Y(v,x_5)w,x_3)\phi \no
\eeqa
where $P_{x_5^{-n-1}}$ projects onto the formal series with only
the terms having the factor $x_5^{-n-1}$. Observe that the lowest
power of $x_5$ is $-n_0-1$ where $n_0$ is the highest integer such
that $v_{n_0}w\neq 0$. Continuing, we find
\beqa
    &&    Y_M(Y(Y(u,x_0)v,x_5)w,x_3)\phi \no \\ && =
    \lim_{x_4\to x_5+x_0}
    \lt(\frc{x_4-x_5}{x_0}\rt)^{k(u,v)}
    \sum_{n\in\Z}\quad
    \lim_{x_1^{\frc1p}\to (x_3+x_4)^{\frc1p}}
    \lt(\frc{x_1-x_3}{x_4}\rt)^{k(u,v_nw)} \cdot \no\\ &&\qquad
    \cdot P_{x_5^{-n-1}}
    \lim_{x_2^{\frc1p}\to(x_3+x_5)^{\frc1p}}
    \lt(\frc{x_2-x_3}{x_5}\rt)^{k(v,w)}
    Y_M(u,x_1)Y_M(v,x_2)Y_M(w,x_3)\phi \no
    \no\\ &&=
    \lim_{x_4\to x_5+x_0}
    \lt(\frc{x_4-x_5}{x_0}\rt)^{k(u,v)}
    \sum_{n\in\Z}\quad
    \lim_{x_1^{\frc1p}\to (x_3+x_4)^{\frc1p}}
    \lt(\frc{x_1-x_3}{x_4}\rt)^{k(u,v_nw)} \cdot \no\\ &&\qquad
    \cdot P_{x_5^{-n-1}}
    \lim_{x_2^{\frc1p}\to(x_3+x_5)^{\frc1p}}
    \lt(\frc{x_2-x_3}{x_5}\rt)^{k(v,w)}
    F(u,v,w,\phi;x_1-x_2,x_1-x_3,x_2-x_3,x_1,x_2,x_3) \no
    \no\\ &&=
    \lim_{x_4\to x_5+x_0}
    \lt(\frc{x_4-x_5}{x_0}\rt)^{k(u,v)}
    \sum_{n\in\Z}\quad
    \lim_{x_1^{\frc1p}\to (x_3+x_4)^{\frc1p}}
    \lt(\frc{x_1-x_3}{x_4}\rt)^{k(u,v_nw)} \cdot \no\\ &&\qquad
    \cdot P_{x_5^{-n-1}}
    F(u,v,w,\phi;x_1-x_3-x_5,x_1-x_3,x_5,x_1,x_3+x_5,x_3)~. \no
\eeqa
Using the particular value of $k(u,v_nw)$ given by (\ref{kdesc}),
we can continue evaluating the limits:
\beqa
    &&    Y_M(Y(Y(u,x_0)v,x_5)w,x_3)\phi \no \\ && =
    \lim_{x_4\to x_5+x_0}
    \lt(\frc{x_4-x_5}{x_0}\rt)^{k(u,v)}
    F(u,v,w,\phi;x_4-x_5,x_4,x_5,x_3+x_4,x_3+x_5,x_3) \no \\
    && = \label{resprwasst}
    F(u,v,w,\phi;x_0,x_5+x_0,x_5,x_3+x_5+x_0,x_3+x_5,x_3)~.
\eeqa

Now, consider the following formal series, and use similar
arguments as those above in order to evaluate it:
\beqa
    &&
    \sum_{n\in\Z} x_0^{-n-1} \lim_{x_2^{1/p}\to
    (x_3+x_5)^{1/p}} \lt(\frc{x_2-x_3}{x_5}\rt)^{k(u_nv,w)}
    Y_M(u_nv,x_2)Y_M(w,x_3)\phi
    \no \\ &=&
    \sum_{n\in\Z} \quad \lim_{x_2^{1/p}\to
    (x_3+x_5)^{1/p}} \lt(\frc{x_2-x_3}{x_5}\rt)^{k(u_nv,w)}
    P_{x_0^{-n-1}} Y_M(Y(u,x_0)v,x_2)Y_M(w,x_3)\phi
    \no \\ &=&
    \sum_{n\in\Z} \quad \lim_{x_2^{1/p}\to
    (x_3+x_5)^{1/p}} \lt(\frc{x_2-x_3}{x_5}\rt)^{k(u_nv,w)} \cdot
    \no\\ &&
    \cdot P_{x_0^{-n-1}}
    \lim_{x_1^{1/p}\to(x_2+x_0)^{1/p}}
    \lt(\frc{x_1-x_2}{x_0}\rt)^{k(u,v)}
    Y_M(u,x_1)Y_M(v,x_2)Y_M(w,x_3)\phi
    \no \\ &=&
    \sum_{n\in\Z} \quad \lim_{x_2^{1/p}\to
    (x_3+x_5)^{1/p}} \lt(\frc{x_2-x_3}{x_5}\rt)^{k(u_nv,w)} \cdot
    \no\\ &&
    \cdot P_{x_0^{-n-1}}
    \lim_{x_1^{1/p}\to(x_2+x_0)^{1/p}}
    \lt(\frc{x_1-x_2}{x_0}\rt)^{k(u,v)}
    F(u,v,w,\phi;x_1-x_2,x_1-x_3,x_2-x_3,x_1,x_2,x_3) \no
    \no \\ &=&
    \sum_{n\in\Z} \quad \lim_{x_2^{1/p}\to
    (x_3+x_5)^{1/p}} \lt(\frc{x_2-x_3}{x_5}\rt)^{k(u_nv,w)} \cdot
    \no\\ &&
    \cdot P_{x_0^{-n-1}}
    F(u,v,w,\phi;x_0,x_2+x_0-x_3,x_2-x_3,x_2+x_0,x_2,x_3) \no
    \no \\ &=&
    F(u,v,w,\phi;x_0,x_5+x_0,x_5,x_3+x_5+x_0,x_3+x_5,x_3) ~.
\eeqa
Comparing with (\ref{resprwasst}), and taking any fixed power of
$x_0$, we have
\beq
    Y_M(Y(u_nv,x_5)w,x_3) = \lim_{x_2^{1/p}\to
    (x_3+x_5)^{1/p}} \lt(\frc{x_2-x_3}{x_5}\rt)^{k(u_nv,w)}
    Y_M(u_nv,x_2)Y_M(w,x_3)\phi~.
\eeq
This is modified weak associativity (\ref{modwasst}) for $s=0$ and
for the pairs $u_nv,w$, for all $n\in\Z$. It is a simple matter to
compare, instead, the expression
\[
    Y_M(Y(\nu^{-s}Y(u,x_0)v,x_5)w,x_3)\phi
\]
with the expression
\[
    \sum_{n\in\Z} x_0^{-n-1} \lim_{x_2^{1/p}\to\om^s
    (x_3+x_5)^{1/p}} \lt(\frc{x_2-x_3}{x_5}\rt)^{k(u_nv,w)}
    Y_M(u_nv,x_2)Y_M(w,x_3)\phi
\]
using similar steps, and using the fact that $\nu$ is an
automorphism of $V$. We obtain modified weak associativity
(\ref{modwasst}) for arbitrary $s\in\Z$ and for the pairs
$u_nv,w$, for all $n\in\Z$.

Repeating the argument with $v$ replaced by $\t{u}_nv$ for all
$\t{u}\in U$ and for all $n\in\Z$, and so on, we find modified
weak associativity for all pairs $v,w$ with $v\in V$ and $w\in V$.
This proves the last part of the theorem. \eproof

\begin{remark}
As usual, it is possible to specialize the theorems above to the
case $p=1$ in order to obtain theorems applying to untwisted modules.
\end{remark}

\begin{remark}
In the theorem above (in close relation with theorems of \cite{Li2}), concerned with twisted modules for vertex
operators algebras, we assumed the existence of a vertex operator
algebra $V$. With slight adjustments, the theorem can be made
into a construction theorem for vertex operator algebras themselves,
in relation with constructions of \cite{Li1} (a more complete construction theory,
the representation theory of vertex operator algebras of \cite{Li1}, is explained at length in \cite{LL},
cf. Theorems 5.7.6, 5.7.11 for instance).

\end{remark}

\sect{Proof of the twisted module structure for $S[\nu]$ and
construction of the twisted vertex operator map}
\label{sect-prooftmod}

The twisted module structure of $S[\nu]$, for the vertex operator
$S$, was established in \cite{L1}, \cite{FLM1}, \cite{FLM2} and
\cite{DL2} (assuming that $\nu$ preserves a rational lattice in
$\a{h}$). In this section, we present a new proof of the twisted
module structure of $S[\nu]$ (which does not require this minor
assumption), and we construct explicitly the twisted vertex
operators using modified weak associativity (and in particular, we
calculate the explicit form of $\Delta_x$ in a simple way).

Starting from (\ref{tvertexalphax}), we will construct the twisted
vertex operator $Y_{S[\nu]}(u,x)$ for all $u\in S$. As stated
above, it is usual in studying vertex operator algebras to
construct vertex operators associated to elements of the vertex
operator algebra from vertex operator associated to ``simpler''
elements, by using associativity. However, for twisted vertex
operators, the natural weak associativity (\ref{wasst}), that is
immediately obtained from the Jacobi identity, is somewhat
complicated by the projection operator and hides the simple
structure of the construction. On the other hand, the modified
weak associativity (\ref{modwasst}) is simpler, especially when
written in the form (specialized for convenience to $s=0$)
\beq
    Y_M(Y(u,x_0)v,x_2) = \lim_{x_1^{1/p} \to (x_2+x_0)^{1/p}}
    \lt( \lt(\frc{x_1-x_2}{x_0}\rt)^{k(u,v)} Y_M(u,x_1) Y_M(v,x_2)
    \rt)~.
    \label{modwasst2}
\eeq
In order to better understand this formula, note that, as remarked
in Remark \ref{remares}, the pre-factor
$\lt(\frc{x_1-x_2}{x_0}\rt)^{k(u,v)}$ gives exactly $1$ when the
limit procedure $\lim_{x_1^{1/p} \to (x_2+x_0)^{1/p}}$ is applied
on it alone. The formula above cannot be simplified by replacing
this pre-factor by 1, however, because the limit is not well
defined on the other factor $Y_M(u,x_1) Y_M(v,x_2)$ alone. The
construction principle will be to replace this product of vertex
operators by a normal-ordered product plus extra terms. On the
normal-ordered product, the limit is well defined, so that the
pre-factor multiplying the normal-ordered product can be set to 1.

It is immediate to see that the set $U = \{{\bf 1},\alpha(-1){\bf
1}|\alpha\in\a{h}\}$ is a generating subset for the vertex operator
algebra $S$. Hence, in the proof of the twisted module structure
of $S[\nu]$ using Theorems \ref{theoconstr1} and \ref{theoequil}, we need first to prove weak commutativity (\ref{wcomt}) 
(and it turns out that it holds with
$k(u,v)=2$ -- the integer $k(u,v)=2$ is the one involved in the
corresponding weak commutativity for the vertex operator algebra
$S$) for the operators defined by (\ref{tvertexalphax}) with
$n=1$. Then, we need to construct explicitly the map $Y_{S[\nu]}(\cdot,x)$ for all
elements of $S$ recursively by using modified weak associativity.
Finally, we need to check modified weak associativity (\ref{modwasst}) for $u\in U$ and $v\in V$.
All other requirements of Theorem \ref{theoequil} are
immediate to see from the construction of Section
\ref{sect-Heisenberg}.

\begin{theorem}\label{theogenelemS}
The operators (\ref{tvertexalphax}) with $n=1$ satisfy weak
commutativity (\ref{wcomt}) with $k(u,v)=2$.
\end{theorem}

\proof We have
\beq
    \alpha^\nu(x_1)\beta^\nu(x_2) =
    \: \alpha^\nu(x_1)\beta^\nu(x_2)\: + h(\alpha,\beta,x_1,x_2)
\eeq
with
\beqa
    h(\alpha,\beta,x_1,x_2) &=& \sum_{m\in\frc1p\Z,\,m>0} x_1^{-m-1}
    x_2^{m-1} m\<\alpha_{(mp)},\beta\> \no\\
        &=&\label{fcth}
    \sum_{r=1}^p \lt(\frc{x_2}{x_1}\rt)^{\frc{r}p}\;
    \frc{1-\frc{r}p + \frc{r}p
    \frc{x_1}{x_2}}{(x_1-x_2)^2}\<\alpha_{(r)},\beta\>~.
\eeqa
It is a simple matter to see that
\beq
    (x_1-x_2)^2h(\alpha,\beta,x_1,x_2)=
    \sum_{r=1}^p \lt(\frc{x_2}{x_1}\rt)^{\frc{r}p}\;
    \lt(1-\frc{r}p + \frc{r}p
    \frc{x_1}{x_2}\rt)\<\alpha_{(r)},\beta\>~.
\eeq
is equal to $(x_1-x_2)^2h(\beta,\alpha,x_2,x_1)$, using
$\<\beta_{(r)},\alpha\> = \<\alpha_{(-r)},\beta\>$. Hence
\[
    (x_1-x_2)^2 [\alpha^\nu(x_1),\beta^\nu(x_2)] =0
\]
which proves the theorem. \eproof

\begin{theorem}\label{theotvom}
The space $S[\nu]$ has the structure of a twisted module for the
vertex operator algebra $S$.
The general form of the twisted vertex operator $Y_{S[\nu]}(u,x)$,
for any element $u$ of the vertex operator algebra $S$, is:
\beqa\label{genformY}
    &&    Y_{S[\nu]}(\alpha_j(-n_j)\cdots \alpha_1(-n_1){\bf 1},x)
    \\ &&
    \qquad = \sum_{J\subset \{1,\ldots,j\}}
    f_{\{1,\ldots,j\}\backslash J}(x) \;
    \: \prod_{l\in J} \frc1{(n_l-1)!} \lt( \frc{d}{dx} \rt)^{n_l-1} \alpha_l^\nu(x) \:
    \no
\eeqa
with $n_1,\ldots,n_j \in \Z_+$, for some factors $f_{I}(x)$, where we just write the dependence on
the index set $I$, but that really depend on the elements
$\alpha_i$ and the integer numbers $n_i$ for all $i\in I$. The set
$J$ on which we sum takes the values $\varnothing$ (the empty
set), $\{1,\ldots,j\}$ (if it is different from $\varnothing$),
and all other proper subsets of $\{1,\ldots,j\}$ (if any). The factors
$f_I(x)$ are given by
\beq\label{solfI}
    f_I(x) = \lt\{\begin{array}{ll} 0 & |I| \ {\rm odd} \\
    \displaystyle
    \sum_{s\in {\rm Pairings}(I)} \prod_{l=1}^{|I|/2}
    g_{s_l}(x) & |I| \ {\rm even} \end{array} \rt.
\eeq
where $|I|$ is the cardinal of $I$, where ${\rm Pairings}(I)$ is the set of all distinct sets
$s=\{s_1,\ldots,s_{|I|/2}\}$ of distinct (without any element in
common) pairs $s_l=(i_l,i_l')$ (where the order of elements is not
important) of elements $i_l\neq i_l'$ of $I$ such that
$\{i_1,\ldots,i_{|I|/2},\,i_1',\ldots,i_{|I|/2}'\} = I$, and where
\beq
    g_{(i,i')}(x) = g(\alpha_i,n_i,\alpha_{i'},n_{i'},x)
\eeq
with
\beq\label{fctg}
    g(\alpha,m,\beta,n,x) = \Res_{x_0} \Res_{x_2}
    x_0^{-m}x_2^{-n} \sum_{r=1}^p
    \lt(\frc{x+x_2}{x+x_0}\rt)^{\frc{r}p}\;
    \frc{1-\frc{r}p + \frc{r}p
    \frc{x+x_0}{x+x_2}}{(x_0-x_2)^2}\<\alpha_{(r)},\beta\>~.
\eeq
\end{theorem}

\proof Clearly, for $j=0$ and $j=1$ the form (\ref{genformY}) is consistent with modified weak associativity and well defined
on $S$, and we must have $f_\varnothing(x) = 1$ and $f_I(x)=0$ if
$|I|=1$. Assume
(\ref{genformY}) to be valid for $j$ replaced with $j-1$. With $k$
a nonnegative integer large enough and $n_1,\ldots,n_j\in\Z_+$, we have
\beqa
    && \label{constr}
        Y_{S[\nu]}(\alpha_j(-n_j)\cdots \alpha_1(-n_1){\bf 1},x)
    \\ && \qquad
    =   \Res_{x_0} x_0^{-n_j} Y_{S[\nu]}(Y(\alpha_j(-1){\bf 1},x_0)\alpha_{j-1}(-n_{j-1})
        \cdots \alpha_1(-n_1){\bf 1},x)
    \no\\ && \qquad
    =   \Res_{x_0} x_0^{-n_j} \lim_{x_1^{1/p} \to (x+x_0)^{1/p}}
        \Bigg( \lt(\frc{x_1-x}{x_0}\rt)^{k} \cdot
    \no\\ && \qquad \qquad\qquad\qquad
        Y_{S[\nu]}(\alpha_j(-1){\bf 1},x_1) Y_{S[\nu]}(\alpha_{j-1}(-n_{j-1})
        \cdots \alpha_1(-n_1){\bf 1},x) \Bigg)
    \no\\ && \qquad
    =   \Res_{x_0} x_0^{-n_j} \lim_{x_1^{1/p} \to (x+x_0)^{1/p}}
        \Bigg( \lt(\frc{x_1-x}{x_0}\rt)^{k}
        \cdot
    \no\\ && \qquad \qquad
        \alpha_j^\nu(x_1) \sum_{J\subset \{1,\ldots,j-1\}}
        f_{\{1,\ldots,j-1\}\backslash J}(x) \;
        \: \prod_{l\in J} \frc1{(n_l-1)!} \lt( \frc{d}{dx} \rt)^{n_l-1} \alpha^\nu_l(x) \:
        \Bigg)~.
        \no
\eeqa
Now, using the commutation relations (\ref{commutrel}), it is a
simple matter to obtain
\beqa
    &&
    \alpha^\nu_j(x_1) \: \prod_{l\in J} \lt( \frc{d}{dx} \rt)^{n_l-1} \alpha^\nu_l(x) \:
    \no\\ && \qquad
    = \: \alpha^\nu_j(x_1) \prod_{l\in J} \lt( \frc{d}{dx} \rt)^{n_l-1} \alpha^\nu_l(x) \:
    \no \\ && \qquad\quad
    + \sum_{i\in J} \lt( \frc{\d}{\d x} \rt)^{n_i-1}h(\alpha_j,\alpha_i,x_1,x)\;
        \: \prod_{l\in J\backslash \{i\}} \lt( \frc{d}{dx} \rt)^{n_l-1} \alpha^\nu_l(x) \:
        \no
\eeqa
with $h(\alpha,\beta,x_1,x_2)$ defined in (\ref{fcth}). Note that,
comparing with (\ref{fctg}), we have
\beq
    g(\alpha,m,\beta,n,x) = \Res_{x_0}
    x_0^{-m} \lim_{x_1^{1/p} \to (x+x_0)^{1/p}}
        \Bigg( \lt(\frc{x_1-x}{x_0}\rt)^{k}\Res_{x_2} x_2^{-n}
        h(\alpha,\beta,x_1,x+x_2) \Bigg)~.
\eeq
Using the relation
\beq
    \frc1{(n-1)!} \lt(\frc{d}{dx} \rt)^{n-1}f(x) =
    \Res_{x_2} x_2^{-n} f(x+x_2)
\eeq
for formal series $f(x)$ with finitely many negative powers of
$x$ and for $n\in\Z_+$, we can now evaluate the limit and the residue on the
right-hand side of the last equality of (\ref{constr}):
\beqa\
    && \label{constr2}
        Y_{S[\nu]}(\alpha_j(-n_j)\cdots \alpha_1(-n_1){\bf 1},x)
    \\ && \qquad
    =
    \sum_{J\subset \{1,\ldots,j-1\}}f_{\{1,\ldots,j-1\}\backslash
    J}(x)\cdot
    \no \\ && \qquad\qquad
    \Bigg(\: \frc1{(n_j-1)!} \lt( \frc{d}{dx} \rt)^{n_j-1} \alpha^\nu_j(x)
        \prod_{l\in J} \frc1{(n_l-1)!} \lt( \frc{d}{dx} \rt)^{n_l-1} \alpha^\nu_l(x) \:
    \no\\ && \qquad\qquad
    + \sum_{i\in J} g(\alpha_j,n_j,\alpha_i,n_i,x)\;
        \: \prod_{l\in J\backslash \{i\}} \frc1{(n_l-1)!}
            \lt( \frc{d}{dx} \rt)^{n_l-1} \alpha^\nu_l(x) \:
            \Bigg)~. \no
\eeqa
This is still of the form (\ref{genformY}). Moreover, it is simple to
understand, from comparing (\ref{constr2}) with (\ref{constr}),
that the solution (\ref{solfI}) is correct. 

This construction certainly gives us a map $Y_{S[\nu]}(\cdot,x)$ on the vector space $S\simeq S(\hat{\a{h}}^-)$. We need to verify that this map satisfies
modified weak associativity (\ref{modwasst}). This requires three steps.

First, we need to check that operators on the right-hand side of (\ref{genformY}) are independent of the order of the pairs
$(\alpha_1,n_1),\ldots,(\alpha_j,n_j)$ for any positive integers $n_1,\ldots,n_j$ and any elements $\alpha_1,\ldots,\alpha_j$ of $\a{h}$.
The sum over pairings has this symmetry, and we need to check that
\beq\label{symg}
    g(\alpha,m,\beta,n,x) = g(\beta,n,\alpha,m,x).
\eeq
This is not immediately obvious, because, by the binomial
expansion convention, $(x_0-x_2)^{-2} \neq (x_2-x_0)^{-2}$).
This symmetry can indeed be checked:
\beqa
    &&
    g(\alpha,m,\beta,n,x) - g(\beta,n,\alpha,m,x)
    \\ && \qquad
    = \Res_{x_0} \Res_{x_2}
    x_0^{-m}x_2^{-n} \sum_{r=1}^p\cdot
    \no \\ && \qquad\qquad\cdot
    \lt(\frc{x+x_2}{x+x_0}\rt)^{\frc{r}p}\;
    \lt(1-\frc{r}p + \frc{r}p
    \frc{x+x_0}{x+x_2}\rt)
    \lt(\frc1{(x_0-x_2)^2}-\frc1{(x_2-x_0)^2}\rt)\<\alpha_{(r)},\beta\>
    \no \\ && \qquad
    = \Res_{x_0} \Res_{x_2}
    x_0^{-m}x_2^{-n} \sum_{r=1}^p\cdot
    \no \\ && \qquad\qquad\cdot
    \lt(\frc{x+x_2}{x+x_0}\rt)^{\frc{r}p}\;
    \lt(1-\frc{r}p + \frc{r}p
    \frc{x+x_0}{x+x_2}\rt)
    x_0^{-1} \frc{\d}{\d x_2} \delta\lt(\frc{x_2}{x_0}\rt)\<\alpha_{(r)},\beta\>
    \no \\ && \qquad =0
\eeqa
where in the last step, we moved the derivative $\frc{\d}{\d x_2}$
towards the left using Leibniz's rule, and we used the formal
delta-function property and the fact that $m,n\in\Z_+$.

Second, we need to check that modified weak associativity (\ref{modwasst}) with
$Y_{S[\nu]}(Y(\alpha_j(-1) {\bf 1},x_0) \alpha_{j-1}(-n_{j-1}) \cdots \alpha_1(-n_1) {\bf 1},x)$ for its left-hand side is in agreement, at negative
powers of $x_0$, with
\beqa &&
	Y_{S[\nu]}(\alpha_j(n_j) \alpha_{j-1}(-n_{j-1}) \cdots \alpha_1(-n_1) {\bf 1},x) = \\ &&
	\qquad \sum_{i=1}^{j-1} n_j \delta_{n_j,n_i} \langle\alpha_j,\alpha_i\rangle Y_{S[\nu]}(\alpha_{j-1}(-n_{j-1})
	\cdots \widehat{\alpha_i(-n_i)} \cdots \alpha_1(-n_1) {\bf 1},x) \no
\eeqa
for $n_j\in \N$ (the nonnegative integers)
where $\widehat{\alpha_i(-n_i)}$ means that the operator $\alpha_i(-n_i)$ is omitted. Repeating the derivation (\ref{constr}), (\ref{constr2}),
we see that this is equivalent to requiring (with the definition (\ref{fctg}))
\beq
	g(\alpha,-m,\beta,n) = m\delta_{m,n} \langle \alpha,\beta\rangle
\eeq
for $m\in\N, n\in\Z_+$. This is a consequence of the fact that
\beq\label{idd}
	\lt(\frc{x+x_2}{x+x_0}\rt)^s \frc{1-s+s\frc{x+x_0}{x+x_2}}{(x_0-x_2)^2} = \frc1{(x_0-x_2)^2}+ \C [[x_0,x_2,x^{-1}]]
\eeq
for any $s\in\C$. The quantity
\[
	\lt(\frc{x+x_2}{x+x_0}\rt)^s \frc{1-s+s\frc{x+x_0}{x+x_2}}{(x_0-x_2)^2} - \frc1{(x_0-x_2)^2}
\]
obviously has only nonnegative powers of $x_2$ and nonpositive powers of $x$. On the other hand, it is equal to
\beqa &&
	\lt(\frc{x+x_2}{x+x_0}\rt)^s \frc{1-s+s\frc{x+x_0}{x+x_2}}{(x_2-x_0)^2} - \frc1{(x_2-x_0)^2} + \no \\ && \qquad
	\lt( \lt(\frc{x+x_2}{x+x_0}\rt)^s \lt(1-s+s\frc{x+x_0}{x+x_2}\rt) - 1\rt) x_0^{-1} \frc{\d}{\d x_2}\delta\lt(\frc{x_2}{x_0}\rt) . \no
\eeqa
The first two terms obviously have only nonnegative powers of $x_0$. The third term can be evaluated using Leibniz's rule and gives zero.
This completes the proof of (\ref{idd}).

Third, the structure of the construction of Section \ref{sect-Heisenberg} shows that (\ref{modwasst}) is valid as well for $s\neq0$ (equation (\ref{transnu}) is satisfied).

The other requirements of Theorems \ref{theoequil} can be checked from the construction of Section \ref{sect-Heisenberg}, and with
Theorems \ref{theoconstr1} and \ref{theogenelemS}, this completes the proof. \eproof

Finally, let us mention that formula (\ref{genformY}) with
(\ref{solfI}) immediately leads to the following formula for the
operator $\Delta_x$ introduced in (\ref{tvertex}):
\beq\label{Deltax}
    \Delta_x = \sum_{q_1,q_2=1}^d \sum_{m,n\in\Z_+} \frc{\b\alpha_{q_1}(m)\b\alpha_{q_2}(n)}{m n}
        g(\b\alpha_{q_1},m,\b\alpha_{q_2},n,x)
\eeq
where we recall that $\b\alpha_q,\,q=1,\ldots,d$ form an orthonormal
basis of $\a{h}$. This was first constructed, in a different form, in \cite{FLM1} and \cite{FLM2}.

\sect{The Lie algebra $\h{\D}^+$} \label{sect-D}

We will now apply modified weak associativity and the results of
the previous section, that $S[\nu]$, constructed in Section
\ref{sect-Heisenberg}, is a twisted module for the vertex operator
algebra $S$, in order to study a certain infinite-dimensional Lie
algebra $\h{\D}^+$ and its representations. This follows closely
the results of \cite{DLMi1} and \cite{DLMi2}, and does not give
new results with respect to these works.

Let $\mathcal{D}$ be the Lie algebra of formal differential
operators on $\mathbb{C}^\times$ spanned by $t^n D^r$, where
$D=t\,\frac{d}{dt}$ and $n \in \mathbb{Z}$, $r \in \N$ (the
nonnegative integers). This Lie algebra has an essentially unique
one-dimensional central extension $\hat{\mathcal{D}}=\mathbb{C}c
\oplus \mathcal{D}$ (denoted in the physics literature by
$\mathcal{W}_{1+\infty}$).

The representation theory of the highest weight modules of
$\hat{\mathcal{D}}$ was initiated in \cite{KR}, where, among other
things, the complete classification problem of the so-called {\em
quasi-finite} representations\footnote{These are representations
with finite-dimensional homogeneous subspaces.} was settled.  The
detailed study of the representation theory of certain subalgebras
of $\hat{\mathcal{D}}$ having properties related to those of
certain infinite--rank ``classical'' Lie algebras was initiated in
\cite{KWY} along the lines of \cite{KR}. In \cite{Bl} and
\cite{M2}, related Lie algebras (and superalgebras) are considered
{}from different viewpoints. As in \cite{DLMi1,DLMi2}, we will
follow these lines and concentrate on the Lie subalgebra
$\hat{\mathcal{D}}^+$ described in \cite{Bl} and recalled below.

View the elements $t^n D^r \;(n \in \mathbb{Z},\; r \in \N)$ as
generators of the central extension $\hat{\mathcal{D}}$. They can
be taken to satisfy the following commutation relations (cf.
\cite{KR}):
\beqa \label{dcom}
    && [t^m f(D), t^n g(D)]= \nn &&
    t^{m+n}(f(D+n)g(D)-g(D+m)f(D))+\Psi(t^m f(D),t^ng(D)) c, \no
\eeqa
where $f$ and $g$ are polynomials and $\Psi$ is the $2$--cocycle
(cf. \cite{KR}) determined by
$$
    \Psi(t^m f(D), t^n g(D))=-\Psi(t^n g(D), t^m f(D))=
    \delta_{m+n,0} \sum_{i=1}^m f(-i)g(m-i), \; m > 0.
$$
We consider the Lie subalgebra $\D^+$ of $\mathcal{D}$ generated
by the formal differential operators
\beq
\label{Lnrdef}
    L_{n}^{(r)}=(-1)^{r+1} D^r (t^n D) D^r ,
\eeq
where $n\in\Z,\;r \in\N$ \cite{Bl}. The subalgebra $\D^+$ has an
essentially unique central extension (cf. \cite{N}) and this
extension may be obtained by restriction of the $2$--cocycle
$\Psi$ to $\mathcal{D}^+$. Let $\hat{\mathcal{D}}^+=\mathbb{C}c
\oplus \mathcal{D}^+$ be the nontrivial central extension defined
via the slightly normalized $2$--cocycle $-\frac{1}{2}\Psi$, and
view the elements $L_{n}^{(r)}$ as elements of
$\hat{\mathcal{D}}^+$. This normalization gives, in particular,
the usual Virasoro algebra bracket relations
\beq\label{virasoro}
    [L_m^{(0)},L_n^{(0)}]=(m-n)L_{m+n}^{(0)}+\frac{m^3-m}{12} \delta_{m+n,0}\,c.
\eeq

In \cite{Bl} Bloch discovered that the Lie algebra
$\hat{\mathcal{D}}^+$ can be defined in terms of generators that
lead to a simplification of the central term in the Lie bracket
relations. Oddly enough, if we let
\beq\label{bLnrbloch}
    \bar{L}_n^{(r)}=L_n^{(r)}+ \frac{(-1)^r}{2}\zeta(-1-2r)\delta_{n,0}c,
\eeq
then the central term in the commutator
\beq \label{bl2coc}
    [\bar{L}_{m}^{(r)},\bar{L}_{n}^{(s)}]= \sum_{i={\rm min}(r,s)}^{r+s}
    a_{i}^{(r,s)}(m,n)
    \bar{L}_{m+n}^{(i)}+ \frc{(r+s+1)!^2}{2(2(r+s)+3)!} m^{2(r+s)+3}  \delta_{m+n,0} \, c
\eeq
is a pure monomial (here  $a_{i}^{(r,s)}(m,n)$ are structure
constants), in contrast to the central term in (\ref{virasoro})
and in other bracket relations that can be found {}from
(\ref{dcom}). As was announced in \cite{L3}, \cite{L4} and shown
in \cite{DLMi2}, in order to conceptualize this simplification
(especially the appearance of zeta-values) one can construct
certain infinite-dimensional projective representations of $\D^+$
using vertex operators.

Let us explain Bloch's construction \cite{Bl}. Consider the Lie
algebra $\h{\a{h}}$ introduced in Section \ref{sect-Heisenberg},
and its induced (level-one) module $S$. Then the correspondence
\beq\label{Lnr}
    L_n^{(r)} \mapsto  \frc12 \sum_{q=1}^d \sum_{j\in\Z} j^r (n-j)^r
    \: \b\alpha_q(j) \b\alpha_q(n-j) \: \com{n\in\Z}\,,\;
    c\mapsto d,
\eeq
where we recall that $\{\b\alpha_q\}$ is an orthonormal basis of
$\a{h}$ and $\:\cdot\:$ is the usual normal ordering, gives a
representation of $\hat{\mathcal{D}}^+$. Let us denote the
operator on the right--hand side of (\ref{Lnr}) by $L^{(r)}(n)$.
In particular, the operators $L^{(0)}(m) \;(m\in\Z)$ give a
well-known representation of the Virasoro algebra with central
charge $c\mapsto d$,
$$
    [L^{(0)}(m),L^{(0)}(n)] = (m-n)L^{(0)}(m+n) + d \,\frc{m^3-m}{12}\, \delta_{m+n,0},
$$
and the construction \eq{Lnr} for those operators is the standard
realization of the Virasoro algebra on a module for a Heisenberg
Lie algebra (cf. \cite{FLM2}).

Without going into any detail, let us mention that Bloch \cite{Bl}
also studied certain natural graded traces using this
representation of $\b{L}^{(r)}_0$, and found that, as in the
well-known case of the Virasoro algebra $r=0$, they possess nice
modular properties.

The appearance of zeta--values in \eq{bLnrbloch} can be
conceptualized by the following heuristic argument \cite{Bl}:
Suppose that we remove the normal ordering in (\ref{Lnr}) and use
the relation $[\b\alpha_q(m),\b\alpha_q(-m)]=m$ to rewrite
$\b\alpha_q(m) \b\alpha_q(-m)$, with $m \geq 0$,  as $\b\alpha_q(-m)
\b\alpha_q(m)+m$. It is easy to see that the resulting expression
contains an infinite formal divergent series of the form
$$1^{2r+1}+2^{2r+1}+3^{2r+1}+ \cdots .$$
A heuristic argument of Euler's suggests replacing this formal
expression by $\zeta(-1-2r)$, where $\zeta$ is the (analytically
continued) Riemann $\zeta$--function. The resulting
(zeta--regularized) operator is well defined and gives the action
of $\bar{L}_n^{(r)}$; such operators satisfy the bracket relations
(\ref{bl2coc}).

\ssect{Realization in $S$: zeta function at negative integers}

In order to understand the appearance of the zeta function at
negative integers using the vertex operator algebra $S$, following
\cite{L3,L4,DLMi1,DLMi2}, we need to introduce slightly different
vertex operators. Consider a vertex operator algebra $V$. The {\em
homogeneous vertex operators} are defined by
\beq\label{homo}
    X(u,x) = Y(x^{L(0)}u,x)\quad (u\in V)~.
\eeq
The most important property of these operators, for us, is the
homogenous version of modified weak associativity:
\beq\label{hmodwass}
    \lim_{x_1 \to e^y x_2}
    \lt(\lt(\frc{x_1}{x_2}-1\rt)^{k(u,v)} X(u,x_1)X(v,x_2)\rt) = \lt(e^y-1\rt)^{k(u,v)}
        X(Y[u,y]v,x_2)
\eeq
for $u,v\in V$, and $k(u,v)$ as in Theorem \ref{theomwa}.
Interestingly, in this relation, yet a new type of vertex operator
appears:
\beq\label{cylY}
    Y[u,y] = Y(e^{yL(0)}u,e^y-1)\quad (u\in V)~.
\eeq
This vertex operator map generates a vertex operator algebra that
is {\em isomorphic} to $V$, and geometrically corresponding to a
change to {\em cylindrical coordinates}. Those properties were
proven in \cite{Z1,Z2}.

Consider the vertex operator algebra $S$. Recall that the Virasoro
generators $L(n)$ acting on $S$ are given by the operators on the
right-hand side of (\ref{Lnr}) with $r=0$. It is a simple matter
to verify that $X(\alpha(-1){\bf 1},x) = \alpha\<x\>$, with
\beq
    \alpha\<x\> = \sum_{n\in \Z} \alpha(n) x^{-n}~.
\eeq
Consider now the following formal series, acting on $S$:
\beq \label{iterbar}
    {\bar{L}}^{y_1,y_2}\<x\>= X
    \lt(\frac{1}{2} \sum_{q=1}^d Y[\b\alpha_q(-1) {\bf 1} ,y_1-y_2]
    \b\alpha_q(-1) {\bf 1} ,e^{y_2}x \rt).
\eeq
By (\ref{hmodwass}), we have
\beq\label{Lbarlimit}
    \b{L}^{y_1,y_2}\<x_2\>
    =
    \frc12 \lim_{x_1 \to  x_2}
    \sum_{q=1}^d \lt(\lt( \frac{ \frac{x_1}{x_2}e^{y_1-y_2} -1 }
    {e^{y_1-y_2} -1} \rt)^k
        \b\alpha_q \< e^{y_1}x_1\> \b\alpha_q\< e^{y_2}x_2\>
       \rt)
\eeq
for any fixed $k\in\N,\,k\ge2$. Using
\[
    \sum_{q=1}^{d}
    \b\alpha_q \< e^{y_1}x_1\> \b\alpha_q\<e^{y_2}x_2\>
    = \sum_{q=1}^{d} \:\b\alpha_q \< e^{y_1}x_1\>
    \b\alpha_q\<e^{y_2}x_2\>\: -
    \frc{\d}{\d y_1} \lt(\frc{1}{
    1-\frc{x_2}{x_1}e^{-y_1+y_2}} \rt)~,
\]
we immediately find that
\beq \label{bLy1y2x}
    {\bar{L}}^{y_1,y_2}\<x\> = \frc12 \sum_{q=1}^{d}
    \: \b\alpha_q\langle e^{y_1}x\rangle
    \b\alpha_q\langle e^{y_2}x\rangle  \:  - \frc12 \frc{\d}{\d y_1} \lt(
        \frac{1}{1-e^{-y_1+y_2}} \rt).
\eeq
Defining the operators
$\b{L}^{r_1,r_2}(n),\,r_1,r_2\in\N,\,n\in\Z$ via
\beq\label{bLr1r2n}
    \b{L}^{y_1,y_2}\<x\>= \frc{1}2 \frc{d}{(y_1-y_2)^2} +
        \sum_{n\in\Z,\;r_1,r_2\,\in\N} \b{L}^{r_1,r_2}(n) x^{-n}
        \frc{y_1^{r_1}y_2^{r_2}}{r_1!r_2!}~,
\eeq
it is simple to see, using (\ref{Lnr}), that the correspondence
\beq\label{repbLnr}
    \b{L}_n^{(r)} \mapsto \b{L}^{r,r}(n),\quad c\mapsto d
\eeq
for $n\in\Z,\,r\in\N$ gives a representation of the generators
(\ref{bLnrbloch}) of the algebra $\h{\D}^+$. Recall that these
generators were introduced by Bloch in order to simplify the
central term in the commutation relations.

As was shown in \cite{DLMi2}, from the expression
(\ref{Lbarlimit}) of the formal series $L^{y_1,y_2}\<x\>$,
involving formal limits, it is a simple matter to compute the
following commutators, first written in \cite{L3}:
\begin{eqnarray}\label{Lbarbracketsalg}
    \lefteqn{[{\bar L}^{y_1,y_2}\<x_1\>,{\bar
    L}^{y_3,y_4}\<x_2\>]} \\
    &&= - {\frac{1}{2}} \frac{\partial}{\partial y_1} \biggl({\bar
    L}^{-y_1+y_2+y_3,y_4}\<x_2\> \delta
    \left({\frac{e^{y_1}x_1}{e^{y_3}x_2}}\right)+ {\bar
    L}^{-y_1+y_2+y_4,y_3}\<x_2\>
    \delta \left({\frac{e^{y_1}x_1}{e^{y_4}x_2}}\right)\biggr)\nonumber\\
    &&\quad - {\frac{1}{2}} \frac{\partial}{\partial y_2} \biggl({\bar
    L}^{y_1-y_2+y_3,y_4}\<x_2\> \delta
    \left({\frac{e^{y_2}x_1}{e^{y_3}x_2}}\right)+ {\bar
    L}^{y_1-y_2+y_4,y_3}\<x_2\> \delta
    \left({\frac{e^{y_2}x_1}{e^{y_4}x_2}}\right)\biggr)~. \no
\end{eqnarray}
As was announced in \cite{L3,L4} and explained in \cite{DLMi2}, a
simple analysis of this commutator shows that the central term in
the commutators of the generators (\ref{repbLnr}) is a pure
monomial, as in (\ref{bl2coc}). This gives a simple explanation of
Bloch's phenomenon using the vertex operator algebras $S$.
Moreover, the definition (\ref{iterbar}) says that the operators
$L^{y_1,y_2}(x)$ represent on $V$ the image of some fundamental
algebra elements of $V$ under the transformation to the cylinder.
These fundamental elements, being closely related to the Virasoro
element $\omega$, can be expected, when transformed to the
cylinder, to lead to graded traces with simple modular properties,
in agreement with the observations of Bloch \cite{Bl}.

\ssect{Representations on $S[\nu]$: Bernoulli polynomials at
rational values}

Following \cite{DLMi1,DLMi2}, we will now construct a
representation of $\h{\D}^+$ on $S[\nu]$. The property that a
twisted module is a true module on the fixed-point subalgebra will
be essential below in this construction. This property is
guaranteed by the operator $\Delta_x$ that we calculated above
(\ref{Deltax}). In order to have the correction terms for the
representation of the algebra $\h{D}^+$ on the twisted space
$S[\nu]$, one could apply $e^{\Delta_x}$ on the vectors generating
the representation of the whole algebra $\h\D^+$. This can be a
complicated problem, mainly because generators of
$\hat{\mathcal{D}}^+$ have arbitrarily large weights. In line with
\cite{DLMi1,DLMi2}, below we will calculate the correction terms
directly using the modified weak associativity relation for
twisted operators, as well as the simple result
(\ref{tvertexalphax}). Hence in this argument, the explicit action
of $\Delta_x$ on vectors generating the representation of the
algebra $\h\D^+$ is not of importance; all we need to know is that
{\em there exists} such an operator $\Delta_x$ giving to the space
$S[\nu]$ the properties of a twisted module for the vertex
operator algebra $S$.

In parallel to the previous sub-section, we need to introduce {\em
homogeneous twisted vertex operators}. Being given a vertex
operator algebra $V$ and a $\nu$-twisted $V$-module $M$, they are
defined by
\beq\label{homot}
    X_M(u,x) = Y_M(x^{L(0)}u,x)\quad (u\in V)~.
\eeq
Again, the most important property of these operators, for us, is
the homogenous version of modified weak associativity for twisted
vertex operators:
\beqa \label{hmodwasst}
&&  \lim_{x_1^{1/p} \to \omega_p^s (e^y x_2)^{1/p}}
\lt(\lt(\frc{x_1}{x_2}-1\rt)^{k(u,v)} X_M(u,x_1)X_M(v,x_2)\rt) \nn
&&\qquad = \lt(e^y-1\rt)^{k(u,v)} X_M(Y[\nu^{-s}u,y]v,x_2).
\eeqa
for $u,v\in V,\;s\in\Z$, and $k(u,v)$ as in Theorem
\ref{theomwat}. Recall the definition of $Y[u,y]$ in (\ref{cylY}).

Consider the vertex operator algebra $S$ and its twisted module
$S[\nu]$. It is a simple matter to verify that
$X_{S[\nu]}(\alpha(-1){\bf 1},x) = \alpha^\nu\<x\>$, with
\beq
    \alpha^\nu\<x\> = \sum_{n\in \frc1p\Z} \alpha^\nu(n) x^{-n}~.
\eeq
Consider now the following formal series, acting on $S[\nu]$:
\beq \label{iterbart}
    {\bar{L}}^{\nu;y_1,y_2}\<x\>= X_{S[\nu]}
    \lt(\frac{1}{2} \sum_{q=1}^d Y[\b\alpha_q(-1) {\bf 1} ,y_1-y_2]
    \b\alpha_q(-1) {\bf 1} ,e^{y_2}x \rt).
\eeq
Since the operator $\frac{1}{2} \sum_{q=1}^d Y[\b\alpha_q(-1) {\bf
1} ,y_1-y_2] \b\alpha_q(-1) {\bf 1} $ is in the fixed point
subalgebra of $S$, it is immediate that these operators satisfy
the same commutation relations as (\ref{Lbarbracketsalg}):
\begin{eqnarray}\label{Lbarbracketstmod}
    \lefteqn{[{\bar L}^{\nu;y_1,y_2}\<x_1\>,{\bar
    L}^{\nu;y_3,y_4}\<x_2\>]} \\
    &&= - {\frac{1}{2}} \frac{\partial}{\partial y_1} \biggl({\bar
    L}^{\nu;-y_1+y_2+y_3,y_4}\<x_2\> \delta
    \left({\frac{e^{y_1}x_1}{e^{y_3}x_2}}\right)+ {\bar
    L}^{\nu;-y_1+y_2+y_4,y_3}\<x_2\>
    \delta \left({\frac{e^{y_1}x_1}{e^{y_4}x_2}}\right)\biggr)\nonumber\\
    &&\quad - {\frac{1}{2}} \frac{\partial}{\partial y_2} \biggl({\bar
    L}^{\nu;y_1-y_2+y_3,y_4}\<x_2\> \delta
    \left({\frac{e^{y_2}x_1}{e^{y_3}x_2}}\right)+ {\bar
    L}^{\nu;y_1-y_2+y_4,y_3}\<x_2\> \delta
    \left({\frac{e^{y_2}x_1}{e^{y_4}x_2}}\right)\biggr)~. \no
\end{eqnarray}
Hence, defining the operators
$\b{L}^{\nu;r_1,r_2}(n),\,r_1,r_2\in\N,\,n\in\Z$ via
\beq\label{bLnur1r2n}
    \b{L}^{\nu;y_1,y_2}\<x\>= \frc{1}2 \frc{d}{(y_1-y_2)^2} +
        \sum_{n\in\Z,\;r_1,r_2\,\in\N} \b{L}^{\nu;r_1,r_2}(n) x^{-n}
        \frc{y_1^{r_1}y_2^{r_2}}{r_1!r_2!}
\eeq
(in particular, only integer powers of $x$ appear in this
expansion), we conclude that they satisfy the same commutation
relations as the operators $\b{L}^{r_1,r_2}(n)$ introduced in
(\ref{bLr1r2n}). Then, as in (\ref{repbLnr}), we can expect that
the correspondence
\beq
    \b{L}_n^{(r)} \mapsto \b{L}^{\nu;r,r}(n),\quad c\mapsto d
\eeq
for $n\in\Z,\,r\in\N$ gives a representation of the generators
(\ref{bLnrbloch}) of the algebra $\h{\D}^+$ on $S[\nu]$. Bringing
this expectation to a proof needs a little more analysis (in
particular, one needs to show that the $L^{\nu;r_1,r_2}(n)$ are
related to the $L^{\nu;r,r}(n)$ in the same way as the
$L^{r_1,r_2}(n)$ are related to the $L^{r,r}(n)$), which is done
in detail in \cite{DLMi2}.

Now, by (\ref{hmodwasst}) we have
\beq\label{Lbarlimitnu}
    \b{L}^{\nu;y_1,y_2}\<x_2\>
    =
    \frc12 \lim_{x_1 \to  x_2}
    \sum_{q=1}^d \lt(\lt( \frac{ \frac{x_1}{x_2}e^{y_1-y_2} -1 }
    {e^{y_1-y_2} -1} \rt)^k
        \b\alpha^\nu_q \< e^{y_1}x_1\> \b\alpha^\nu_q\< e^{y_2}x_2\> \rt)
\eeq
for any fixed $k\in\N,\,k\ge2$, which gives
\beq \label{bLnuy1y2x}
    {\bar{L}}^{\nu;y_1,y_2}\<x\>= \frc12 \sum_{q=1}^{d}
    \: \b\alpha^\nu_q\langle e^{y_1}x\rangle
    \b\alpha^\nu_q\langle e^{y_2}x\rangle
    - \frc12 \frc{\d}{\d y_1} \lt(\sum_{k=0}^{p-1}
        \frac{e^{\frac{k(-y_1+y_2)}{p}}
    {\rm dim} \ \goth{h}_{(k)}}
    {1-e^{-y_1+y_2}} \rt)
\eeq
using
\[
    \sum_{q=1}^{d}
    \b\alpha^\nu_q \< e^{y_1}x_1\> \b\alpha^\nu_q\<e^{y_2}x_2\>
    = \sum_{q=1}^{d} \:\b\alpha^\nu_q \< e^{y_1}x_1\>
    \b\alpha^\nu_q\<e^{y_2}x_2\>\: -
    \frc{\d}{\d y_1} \lt(
    \sum_{k=0}^{p-1} \frc{e^{\frc{k(-y_1+y_2)}p} \; {\rm dim} \
    \goth{h}_{(k)}}{
    1-\frc{x_2}{x_1}e^{-y_1+y_2}} \rt).
\]
Evaluating the operators $\b{L}^{\nu;r,r}(n)$ from
(\ref{bLnur1r2n}), we conclude that the operators
\beqa
    \b{L}^{\nu;r,r}(n) &=& \frc12 \sum_{q=1}^d \sum_{j\in \frac{1}{p}
    \Z} j^{r}(n-j)^{r}\,
    \: \b\alpha_q^\nu (j)\b\alpha_q^\nu (n-j) \: \no\\   &&  -
    \delta_{n,0} \, \frc{(-1)^r}{4(r+1)} \sum_{k=0}^{p-1}
    \dim\a{h}_{(k)} B_{2(r+1)}(k/p)
\eeqa
form a representation, on $S[\nu]$, of the generators
(\ref{bLnrbloch}) for the Lie algebra $\h{\D}^+$. Notice the
appearance of the Bernoulli polynomials. From our construction,
this is seen to be directly related to general properties of
homogeneous twisted vertex operators.

The next result is a simple consequence of the discussion above.
It was shown in \cite{DLMi2}. It describes the action of the
``Cartan subalgebra'' of $\hat{\mathcal{D}}^+$ on a highest weight
vector of a canonical quasi-finite $\hat{\mathcal{D}}^+$--module;
here we are using the terminology of \cite{KR}. This corollary
gives the ``correction'' terms referred to in the introduction.

\begin{corol}
Given a highest weight $\hat{\mathcal{D}}^+$--module $W$, let
$\delta$ be the linear functional on the ``Cartan subalgebra'' of
$\hat{\mathcal{D}}^+$ (spanned by $L^{(k)}_{0}$ for $k\in \N$)
defined by
$$
    {L}^{(k)}_{0} \cdot w =(-1)^k\delta\lt( {L}^{(k)}_{0} \rt)w,
$$
where $w$ is a generating highest weight vector of $W$, and let
$\Delta(x)$ be the generating function
$$
    \Delta(x)=\sum_{k \ge 1} \frac{\delta({L}^{(k)}_0)x^{2k}}{(2k)!}
$$
(cf. \cite{KR}). Then for every automorphism $\nu$ of period $p$
as above,
$$
    \mathcal{U}(\hat{\mathcal{D}}^+) \cdot 1 \subset S[\nu]
$$
is a quasi--finite highest weight $\hat{\mathcal{D}}^+$--module
satisfying
\beq
    \Delta(x) = \frac{1}{2} \frac{d}{dx} \sum_{k=0}^{p-1} \frac{
    (e^{\frac{kx}{p}} -1) {\rm dim} \ \goth{h}_{(k)}} {1-e^x}.
\eeq
\end{corol}

\small{

}

\noindent {\small \sc 
Rudolf Peierls Centre for Theoretical
Physics,\\ Oxford University, UK} \\
{\em Current address}:\\
\noindent {\small \sc 
Department of Mathematical Sciences\\
Durham University, UK}\\
{\em E-mail address}: benjamin.doyon@durham.ac.uk

\end{document}